\documentclass[12pt,a4paper]{article}
\usepackage[utf8]{inputenc}
\usepackage[T1]{fontenc}
\usepackage{lmodern}
\usepackage{amsmath,amssymb,amsthm}
\usepackage{stmaryrd}
\usepackage{natbib}
\usepackage{setspace}
\usepackage{geometry}
\usepackage{url}
\usepackage{microtype}

\newcommand{\N}{\mathbb{N}}
\newcommand{\Z}{\mathbb{Z}}

\newcommand{\R}{\mathbb{R}}

\newcommand{\fS}{\mathfrak{S} }
\newcommand{\fB}{\mathfrak{B}}
\newcommand{\fA}{\mathfrak{A}}

\newcommand{\NN}{\mathfrak{N}}
\newcommand{\RR}{\mathfrak{R}}
\newcommand{\M}{\mathfrak{M}}
\newcommand{\Nx}{\mathfrak{N}^\times}

\newcommand{\id}{\mathrm{id}}
\newcommand{\Iso}{\mathrm{Iso}}

\newcommand{\Th}{\mathrm{Th}}

\newcommand{\Lom}{L_{\omega_1^{\mathrm{CK}}}}
\newcommand{\omCK}{\omega_1^{\mathrm{CK}}}

\newcommand{\supp}{\mathrm{supp}}

\newcommand{\then}{\rightarrow}
\newcommand{\gdw}{\leftrightarrow}

\newcommand{\Gdw}{\Leftrightarrow}

\DeclareMathOperator{\PA}{PA}
\DeclareMathOperator{\QA}{Q}
\DeclareMathOperator{\KP}{KP}
\DeclareMathOperator{\Ze}{Z}
\DeclareMathOperator{\ZF}{ZF}
\DeclareMathOperator{\ZFC}{ZFC}
\DeclareMathOperator{\PrA}{PrA}
\DeclareMathOperator{\SucA}{SucA}
\DeclareMathOperator{\RCF}{RCF}
\DeclareMathOperator{\ACF}{\mathrm{ACF}}

\DeclareMathOperator{\SkA}{SkA}
\DeclareMathOperator{\BA}{BA}
\DeclareMathOperator{\F}{F}
\DeclareMathOperator{\Span}{span}

\theoremstyle{plain}
\newtheorem{theorem}{Theorem}[section]
\newtheorem{definition}[theorem]{Definition}

\newtheorem{fact}[theorem]{Fact}

\newtheorem{lemma}[theorem]{Lemma}
\newtheorem{corollary}[theorem]{Corollary}

\title{Rigid Theories}
\author{%
  Mirko Engler\thanks{%
    Department of Philosophy, University of Vienna, Austria.
    Emails: \texttt{mirko.engler@univie.ac.at} (Mirko Engler);
    \texttt{benjamin.zayton@univie.ac.at} (Benjamin Zayton).
  }
  \and
  Benjamin Zayton\footnotemark[1]
}
\date{}

\begin{document}
\maketitle
{\renewcommand{\thefootnote}{}%
 \footnotetext{\emph{Acknowledgements.} Since we first became interested in studying the absence of non-trivial self-interpretations of theories as a notion in its own right, a number of researchers have generously helped us develop this investigation with their expertise. Albert Visser drew our attention to the work of Zoethout and of Pakhomov and Zapryagaev. During a short research stay of the first author in Warsaw, numerous discussions with Mateusz {\L}e{\l}yk, Piotr Gruza, Konstantinos Papafilippou, and Patryk Szlufik helped to develop conceptual connections and examples. After early versions of this article had circulated, Ali Enayat joined the discussion and, together with Mateusz, helped to improve and further develop the results, particularly in set theory. We are indebted to all of them for generously sharing their ideas.}}
\begin{abstract}
\noindent
A structure is rigid if it admits no non-trivial automorphisms. We investigate an analogous invariance property for first-order theories, with relative interpretations playing the role of morphisms. A theory $T$ is called \emph{rigid} if and only if every self-interpretation of $T$ is definably isomorphic to the identity interpretation. Extending work of \cite{Pakhomov2020} on Presburger Arithmetic, we establish rigidity for Successor Arithmetic and Skolem Arithmetic (with respect to one-dimensional interpretations), True Arithmetic, and the theories of fields $\RCF$ and $\ACF$. By contrast, $\PA$ and $\ZF$ admit infinitely many self-interpretations, though both remain rigid with respect to (co-)retractions. Our main negative result is that $\ZFC$ has no consistent rigid extension in its own language. We compare rigidity with categoricity, model rigidity, completeness, and solidity, showing it to be independent of the first three, and indicate the relevance of rigidity to questions about semantic determinacy.
\end{abstract}

\section{Introduction}

A classic topic in model theory is the study of the properties of structures invariant under appropriate transformations. More precisely, the automorphisms of a structure $\M$ are taken to provide conceptual insight into the structure of $\M$. Of particular interest is the case where a structure $\M$ is \emph{rigid}, that is, the case where $\M$ admits no non-trivial automorphisms.

\medskip\noindent
While the study of the automorphism groups of structures has received a lot of attention in model theory, analogous questions for theories have received no sustained attention. Passing from structures to theories, the question arises what the appropriate notion of morphism should be. A natural choice for this is the notion of relative interpretations between theories. Implicitly, these were introduced as morphisms between theories in \cite{tarski1953}, but their role as morphisms in categories of theories and interpretations was first explicitly studied in the foundational work of \cite{Visser2006}. 

\medskip\noindent
Moreover, relative interpretations are analogous to homomorphisms of structures in that they preserve the internal structure of theories, because they respect their inferential structure, just as homomorphisms of models respect extensions of predicates: they map theorems to theorems and preserve the logical structure of formulae. The appropriate analogue of automorphisms depends on the identity conditions one imposes on interpretations, but it is clear that, for instance, bi-interpretations from a theory to itself have a good claim to being symmetries of the theory in just the same way in which automorphisms are symmetries of a model-theoretic structure. 

\medskip\noindent
The relative interpretations of a theory in itself form a monoid under composition, and it is the structure of this monoid that determines the degree of self-translation freedom of a theory. We call a theory \emph{rigid} in the corresponding sense if the quotient of this monoid under  homotopy is trivial; homotopy is here meant in the sense of \citet{hodges1993} or \citet{Visser2006}. That is, the quotient is trivial if every self-interpretation is definably isomorphic to the identity interpretation. In this paper, we aim to advance the study of this concept by developing some initial general theory, and, primarily, by examining mathematically relevant theories from arithmetic, algebra, and set theory with respect to this property. 

\medskip\noindent
Rigidity can be understood as being relative to a class of interpretations and an equivalence relation on that class. In our investigation, the main focus will be on parameter-free, possibly multi-dimensional interpretations which preserve identity, identified up to definable isomorphism. To the best of our knowledge, \cite{Pakhomov2020} were the first to investigate the triviality of self-interpretations in this generality. They demonstrated, using our new terminology, that Presburger Arithmetic is rigid, answering a question of Visser. We extend this investigation to other arithmetical theories and show that Successor Arithmetic and Skolem Arithmetic are rigid with respect to one-dimensional interpretations. Moving up, Peano Arithmetic has infinitely many self-interpretations, but is rigid with respect to (co-)retractions. Its most important completion, True Arithmetic, turns out to be rigid. Furthermore, we deduce the rigidity of $\RCF$ and $\ACF_p$ for all characteristics $p$, as well as of $\ACF$, from established model-theoretic results. In contrast, we will show that $\ZFC$ has no consistent rigid extension in its own language.

\medskip\noindent
The study of the rigidity of a theory also has a philosophical motivation: A relative interpretation of a theory in another theory can be seen as a reading of the first theory's axioms in the second theory. An influential tradition, which can be traced back to \cite{hilbert1899},  \cite{carnap1937}, and \cite{ajdukiewicz1934}, holds that the meaning of a formal theory can be determined by the axioms of that theory. If one takes this seriously, then a self-interpretation of $T$ represents an alternative reading of $T$’s own theorems, which is sanctioned in a precise sense by the axioms of $T$. A theory that allows only trivial self-interpretations is a theory whose axioms leave no room for such alternative readings: the axioms unambiguously determine the meaning of the theorems. Rigidity of theories thus provides a precise, purely syntactic criterion for what one might call the \emph{semantic determinacy} of a theory. A more detailed philosophical examination of the rigidity of theories, as well as their relationship to other concepts of semantic determinacy of theories, will be provided in a separate article.

\section{Rigid Theories}
In this paper, we will be concerned with self-interpretations of a theory $T$ into itself, that is, with relative interpretations $f:T\to T$. We officially restrict ourselves to relational languages, understanding function symbols and constants as relations with appropriate axioms. Throughout, fraktur letters denote structures, and the corresponding italic letter denotes the carrier of the structure: $M$ is the carrier of $\M$, and so on. For the standard models $\NN$ and $\Nx$ of arithmetic the carriers are $\N$ and $\N^\times:=\N_{>0}$, respectively. Moreover, if no parameters are mentioned, we understand \emph{definability} in a structure as \emph{definability without parameters} in that structure. 
\noindent
The definition of a many-dimensional relative interpretation between theories $S$ and $T$ goes as follows:
\begin{definition}[Interpretation between Theories]\label{interpretation}
    Let $S$ and $T$ be first-order theories and $d\geq 1$. For each $n$-ary relation symbol $R$ from the language $\mathcal{L}_S$ of $S$, let $f(R)$ be an $n \times d$-ary formula in $\mathcal{L}_T$, and let $\delta_f$ be a $d$-ary $\mathcal{L}_T$-formula. Set $\bar x_j=(x_{jd},\dots,x_{jd+d-1})$. The assignment $f$ induces a function $\mathcal{L}_S\to\mathcal{L}_T$ via the following recursive clauses: 
    \begin{enumerate}
        \item for all $n,m$: $f(x_n=x_m):=\bigwedge_{i=0}^{d-1}(x_{nd+i}=x_{md+i})$,
        \item for all $n$-ary relation symbols $R$ of $\mathcal{L}_S$:
        \begin{equation*}
            f(R(x_{i_1},\ldots, x_{i_n})):=f(R)(\bar x_{i_1},\dots, \bar x_{i_n}),
        \end{equation*}
        \item for all formulae $\psi, \chi$ in $\mathcal{L}_S$: $f(\neg \psi):=\neg f(\psi)$ and $f(\psi\to \chi):=f(\psi)\to f(\chi)$,
        \item for all formulae $\psi$ in $\mathcal{L}_S$ and all $n$: 
        \begin{equation*}
            f(\forall x_n\psi):=\forall \bar x_n(\delta_f (\bar x_n)\to f(\psi)).
        \end{equation*}
       
    \end{enumerate}
    The function $f$ is an interpretation of $S$ in $T$ if, for all sentences $\varphi$ of $\mathcal{L}_S$, 
    \begin{equation*}
        S\vdash \varphi\Rightarrow T\vdash f(\varphi).
    \end{equation*}
\end{definition}
\noindent
If one allows the identity to be taken to an arbitrary formula for which the theory $T$ proves the translation of the schemata for identity, one obtains the notion of an identity-modifying interpretation. Other notable generalizations are interpretations with parameters and piecewise interpretations. We will not take these into account systematically in this paper.

\medskip\noindent
In general, any interpretation $f:S\to T$ has a contravariantly associated semantic dual
$f^{*}:\mathrm{Mod}(T)\to\mathrm{Mod}(S)$, a functor between the
categories of models with elementary embeddings as morphisms: it assigns to each
$\M\vDash T$ a model $\M_f\vDash S$ with $\M_f\vDash\varphi \Leftrightarrow \M\vDash f(\varphi)$
for every $\mathcal{L}_S$-sentence $\varphi$. If the carrier of a model $\M$ is $M$, we will denote the carrier of the induced model $\M_f$ by $M_f$. We will mostly be interested in interpretations up to homotopy, that is, we will typically identify interpretations $f$ and $g$ if there is a formula $\eta(\bar x,\bar y)$ such that $\eta(\bar x,\bar y)$ defines, in each model $\M\vDash T$, an isomorphism $\M_f\to \M_g$. Syntactically, this amounts to the following:
\begin{definition}[Homotopy between Interpretations]\label{homotopy}
Let $f,g:S\to T$ be interpretations with domain-formulae
$\delta_f(\bar x)$ and $\delta_g(\bar y)$, where $|\bar x|=d_f$ and $|\bar y|=d_g$
are their dimensions. An $\mathcal{L}_T$-formula $\eta(\bar x,\bar y)$
\emph{defines a bijection between $\delta_f$ and $\delta_g$ in $T$} if:
\begin{itemize}
\item $T\vdash \forall \bar x\,\bar y\,\bigl(\eta(\bar x,\bar y)\then \delta_f(\bar x)\land \delta_g(\bar y)\bigr)$
\item $T\vdash \forall \bar x\,\bigl(\delta_f(\bar x)\then \exists \bar y\,\eta(\bar x,\bar y)\bigr)$
\item $T\vdash \forall \bar x\,\bar y\,\bar y'\,\bigl(\eta(\bar x,\bar y)\land \eta(\bar x,\bar y')\then \bar y=\bar y'\bigr)$
\item $T\vdash \forall \bar y\,\bigl(\delta_g(\bar y)\then \exists \bar x\,\eta(\bar x,\bar y)\bigr)$
\item $T\vdash \forall \bar x\,\bar x'\,\bar y\,\bigl(\eta(\bar x,\bar y)\land \eta(\bar x',\bar y)\then \bar x=\bar x'\bigr)$
\end{itemize}
We say $f$ and $g$ are \emph{homotopic} (via $\eta$), written
$f\simeq g$ (via $\eta$) or $\eta:f\Rightarrow g$, if the $\mathcal{L}_T$-formula $\eta$ defines a bijection
between $\delta_f$ and $\delta_g$ in $T$, and, for every $n$-ary relation
symbol $R$ of $\mathcal{L}_S$,
\begin{itemize}
\item $T\vdash \forall \bar x_1\dots \bar x_n\,\forall \bar y_1\dots \bar y_n\,
  \Bigl(\textstyle\bigwedge_{i=1}^n \eta(\bar x_i,\bar y_i)\then
  \bigl(f(R)(\bar x_1,\dots, \bar x_n)\gdw g(R)(\bar y_1,\dots, \bar y_n)\bigr)\Bigr)$,
\end{itemize}
where $|\bar x_i|=d_f$ and $|\bar y_i|=d_g$ for each $i$.
\end{definition}
\noindent
For complete theories, one can check whether a formula defines a homotopy using a single model.
\begin{lemma}[Homotopies in Complete Theories]\label{homotopy_completeness}
    Let $f,g:S\to T$ be interpretations and let $\eta(\bar x,\bar y)$ be an $\mathcal{L}_T$-formula. Assume that $T$ is complete and consistent. Then $f$ and $g$ are homotopic via $\eta$ iff there is a model $\M\vDash T$ such that $\eta$ defines an isomorphism from $\M_f$ to $\M_g$ in $\M$. 
\end{lemma}

\noindent
We will call a language \emph{finite} if its underlying signature contains finitely many relation symbols. If the language of the theory $S$ is finite, the conditions in the definition above can clearly be encoded in a single sentence:
\begin{lemma}[Expressibility of Homotopy for Finite Languages]\label{Iso}
    Let $f,g:S\to T$ be interpretations. Assume that the language $\mathcal{L}_S$ of $S$ is finite. Then, for any $(d_f+d_g)$-ary $\mathcal{L}_T$-formula $\eta$, there is an $\mathcal{L}_T$-sentence $\Iso(\eta,f,g)$ such that 
    \begin{equation*}
        T\vdash \Iso(\eta,f,g) \Gdw \eta:f\Rightarrow g.
    \end{equation*}
    More generally, if $\eta(\bar x,\bar y,\bar z)$ has further free variables $\bar z$, the same conjunction yields a formula $\Iso(\eta,f,g)(\bar z)$ such that for every $\M\vDash T$ and $\bar c\in M^{|\bar z|}$, $\M\vDash\Iso(\eta,f,g)(\bar c)$ iff $\eta(\bar x,\bar y,\bar c)$ defines an isomorphism $\M_f\to\M_g$. 
\end{lemma}

With this idea at hand, completeness and a compactness argument together yield:
\begin{lemma}[Model-Theoretic Criterion for Homotopy]\label{homotopy_compactness}
    Let $f,g:S\to T$ be interpretations, with the language $\mathcal{L}_S$ being finite. Assume that for each model $\M\vDash T$ of $T$, there is a definable isomorphism $\M_f\to \M_g$. Then there is a homotopy $\eta:f\Rightarrow g$. 
\end{lemma}
\noindent
Under composition the self-interpretations $\mathrm{Int}(T,T)$ form a monoid, and
$\simeq$ is a congruence on it; denote the quotient by $\mathrm{End}(T):=\mathrm{Int}(T,T)/{\simeq}$. Its invertible elements are exactly the interpretations which are part of a self-bi-interpretation. The central question of this paper is for which
theories $\mathrm{End}(T)$ is trivial.
\begin{definition}[Rigidity of a Theory]
A theory $T$ is rigid if every self-interpretation $f:T\to T$ is homotopic to the identity interpretation $\id_T$.
\end{definition}
\noindent
Of course, one can also consider variants of the notion of rigidity, e.g., by specializing to a sensible sub-class $C$ of all interpretations and a sensible equivalence relation $\mathrel{E}$ on that class.\footnote{A note on terminology: rigidity in model theory typically refers to symmetries in the guise of automorphisms, but using the term rigidity in connection with an endomorphism monoid is not entirely without precedent, see, e.g., \cite{hedrlin1966rigid} or \cite{kamio2026solution}.} To be more precise, the class $C$ should be such that one obtains a category whose objects are theories and whose morphisms are interpretations in $C$. Concretely, this means that $C$ should contain the identity interpretation for each theory and be closed under composition of interpretations. The equivalence relation $\mathrel{E}$ should be a congruence with respect to composition, that is, if $f_1\mathrel{E}f_2$ and $g_1\mathrel{E}g_2$, then $g_1f_1\mathrel{E}g_2f_2$ whenever these compositions are defined. 
Now, one can call a theory $T$ $C,E$-rigid if any self-interpretation $f:T\to T$ with $f\in C$ is $E$-equivalent to $\id_T$. In such a case, we also say that $T$ is rigid with respect to $C$ up to $E$. Sensible classes $C$ can be generated by a choice of feature of an interpretation (relativization, dimension, identity-preservation, parameter-domain) and by their properties (faithfulness, essential surjectivity, (co-)retraction, part of a bi-interpretation, or a synonymy). The premier examples for equivalence relations for us are provable equivalence, sentential equivalence, equivalence up to pointwise isomorphism of the duals, and equivalence up to uniformly definable isomorphism or homotopy. Typically, if $C$ and $E$ are not specified, $C$ should be taken to be the class of interpretations as in Definition~\ref{interpretation} and $E$ should be taken to be the equivalence relation given by homotopy as in Definition~\ref{homotopy}.\footnote{It should be noted here that for all choices of $C,E$ considered in this paper, inconsistent theories are always $C,E$-rigid.} 

\medskip
\noindent
We now record some basic definitions and technical lemmata, starting with an observation essentially due to Visser. 
\begin{lemma}[Finite Axiomatizability from Right to Left]
If $T$ is a rigid theory in a finite language $\mathcal{L}$ and $T$ is interpretable in a finitely axiomatizable $T' \subseteq T$, then $T$ is finitely axiomatizable.
\end{lemma}
\begin{proof}
     Let $T$ be a rigid theory in a finite language $\mathcal{L}$ and assume that $T$ is interpretable in a finitely axiomatizable subtheory $T'\subseteq T$ via an interpretation $f:T\to T'$. Via composition with the inclusion of $T'$ in $T$, $f$ can be seen as a self-interpretation of $T$ and is hence homotopic to $\id_T$ via a homotopy $\eta$. By Lemma~\ref{Iso}, this can be expressed in one sentence $\Iso(\eta, f,\id_T)$ proved by $T$, because the language $\mathcal{L}$ is finite. There is hence a finitely axiomatized subtheory $S\supseteq T'$ of $T$ such that $S\vdash \Iso(\eta,f,\id_T)$. In particular, for any sentence $\varphi$ in the language $\mathcal{L}_T$, $S\vdash f(\varphi)\leftrightarrow \varphi$. Now, if $T\vdash \varphi$, then $S\vdash f(\varphi)$. But then $S\vdash \varphi$. 
\end{proof}
\noindent
There are some easy preservation results for changes in the class $C$ and the equivalence relation $E$: If $C\subseteq C'$ and $E'\subseteq E$ and $T$ is $C',E'$-rigid, then $T$ is $C,E$-rigid. More interesting is the behaviour of rigidity with respect to theory extension.
\begin{lemma}[Upwards Preservation of Rigidity]\label{upwards_preservation_rigidity}
    If $S$ is rigid and $T$ is an extension of $S$ in the same finite language $\mathcal{L}$, then $T$ is rigid if any of the following conditions hold:
    \begin{itemize}
        \item $S$ is finitely axiomatized.
        \item $T$ is finitely axiomatized relative to $S$, that is, there is a sentence $\psi$ in $\mathcal{L}$ such that $S + \psi =T$. 
    \end{itemize}
\end{lemma}
\begin{proof}
    Assume first that $S$ is rigid and finitely axiomatized. Consider an extension $T$ of $S$ and an interpretation $f:T\to T$. There is a $\psi\in T$ such that $\psi\vdash f(S)$ and $\psi\vdash \exists \bar x\delta_f(\bar x)$. Let $d=d_f$ be the dimension of $f$ and let $\bar x_j:=(x_{jd},\dots,x_{jd+d-1})$. We define the $d$-dimensional diagonal-padded identity $\iota_d:S\to S$ with domain-formula $\Delta_d(\bar x_j):=\bigwedge_{i=1}^{d-1}x_{jd}=x_{jd+i}$ and $\iota_d(R)(\bar x_0,\dots,\bar x_{n-1}):=\bigwedge_{j=0}^{n-1}\Delta_d(\bar x_j)\land R(x_0,x_d,x_{2d},\ldots,x_{(n-1)d})$ for each $n$-ary $R$ in $\mathcal{L}$. Define a disjunctive interpretation $g:S\to S$ with $d_g=d$ by
    \begin{itemize}
        \item $\delta_g(\bar x_0):=(\psi \land \delta_f(\bar x_0))\lor (\neg \psi\land \Delta_d(\bar x_0))$ 
        \item $g(R)(\bar x_0, \dots, \bar x_{n-1}):=(\psi\land f(R)(\bar x_0, \dots, \bar x_{n-1}))\lor (\neg \psi\land \iota_d(R)(\bar x_0,\dots, \bar x_{n-1}))$ \\
        for each $n$-ary $R$ in $\mathcal{L}$.\footnote{This is the identity-preserving analogue of Visser's padding in \cite[§2.2]{Visser2019}.}
    \end{itemize}
    By induction on $\varphi$, one shows that
    \[
    S\vdash g(\varphi)(\bar x_0,\dots, \bar x_{n-1})\leftrightarrow\bigl((\psi\land f(\varphi)(\bar x_0,\dots,\bar x_{n-1}))\lor(\neg\psi\land\iota_d(\varphi)(\bar x_0,\dots, \bar x_{n-1}))\bigr)
    \]
    for all $\mathcal{L}$-formulae $\varphi(x_0,\dots,x_{n-1})$. The formula $\theta(\bar x_0,y):=\Delta_d(\bar x_0)\land y=x_0$ defines a homotopy between $\iota_d$ and $\id_S$, so $\iota_d\simeq\id_S$ and hence $S\vdash\iota_d(\varphi)\leftrightarrow\varphi$ for any $\mathcal{L}$-sentence $\varphi$, whereby $S\vdash g(\varphi)\leftrightarrow\bigl((\psi\land f(\varphi))\lor(\neg\psi\land\varphi)\bigr)$.
    Therefore, if $S \vdash \varphi$, then \begin{align*}
     S + \neg \psi &\vdash \neg \psi \land \varphi \text{, so } S + \neg \psi \vdash g(\varphi)\\
     S +\psi &\vdash \psi \land f(\varphi) \text{, so }  S + \psi \vdash g(\varphi) 
    \end{align*} from which it follows that $g:S\to S$ is an interpretation. Therefore, $g$ is homotopic to $\id_S$ via a formula $\eta$. By Lemma~\ref{Iso}, this can be expressed in one sentence $\Iso(\eta, g,\id_S)$ proved by $S$, because the language $\mathcal{L}$ is finite. Since $S\subseteq T$, $T$ likewise proves $\Iso(\eta,g,\id_S)$. Because $T \vdash \psi$, we have $T\vdash \delta_f(\bar x) \gdw \delta_g(\bar x)$ and $T\vdash f(R)(\bar x_0, \dots, \bar x_{n - 1} ) \gdw g(R)( \bar x_0, \dots, \bar x_{n - 1})$  for all $R$ in $\mathcal{L}$, hence $T\vdash \Iso(\eta,g,\id_S)\leftrightarrow \Iso(\eta,f,\id_T)$, and therefore, $T\vdash \Iso(\eta,f,\id_T)$, whence $f$ is homotopic to the identity over $T$ and $T$ is rigid.

\medskip
\noindent
    For the second case, where $T$ is $S + \psi$, repeat the same construction of a disjunctive interpretation with this $\psi$.
  
\end{proof}
\noindent
Here are some other important properties of interpretations we will consider.
\begin{definition}[Faithful Interpretation]
    An interpretation $f:S\to T$ is faithful iff, for all sentences $\varphi\in \mathcal{L}_S$, 
    \begin{equation*}
        S\vdash \varphi \Leftrightarrow T\vdash f(\varphi).
    \end{equation*}
    
\end{definition}
\begin{definition}[Retractions and Co-Retractions]
    An interpretation $f:S\to T$ is an $E$-retraction iff there is an interpretation $g:T\to S$ such that $f\circ g \mathrel{E}\id_T$. It is a co-retraction or a section iff there is an interpretation $g:T\to S$ such that $g\circ f \mathrel{E} \id_S$.
\end{definition}

\medskip\noindent
Mere finite axiomatizability appears to be too weak for upwards preservation of rigidity with respect to these smaller classes of interpretations, but we still get the following results:
\begin{lemma}[Upwards Preservation of Rigidity for Faithful Interpretations]\label{upwards_preservation_faithful}
    If $S$ is rigid with respect to faithful interpretations and $S\subseteq T$ is an extension of $S$ in the same finite language, then $T$ is rigid with respect to faithful interpretations if $T$ is finitely axiomatized relative to $S$. 
\end{lemma}
\begin{proof}
    Let $T=S+\psi$. Let $f:T\to T$ be faithful. Then define a disjunctive interpretation $g:S\to S$ as in the proof of Lemma~\ref{upwards_preservation_rigidity}, such that $S\vdash g(\varphi)\leftrightarrow((\psi\land f(\varphi))\lor(\neg\psi\land \varphi))$ for $\mathcal{L}$-sentences $\varphi$. It remains to show that $g$ is faithful; the remainder of the proof goes as in Lemma~\ref{upwards_preservation_rigidity}. Now, if $S\vdash g(\varphi)$, then $S\vdash (\psi\land f(\varphi))\lor(\neg\psi\land \varphi)$. This yields $S\vdash \psi\to f(\varphi)$, i.e., $T\vdash f(\varphi)$; since $f$ is faithful over $T$, we get $T\vdash\varphi$, that is, $S\vdash\psi\to\varphi$. The same disjunction also gives $S\vdash\neg\psi\to\varphi$, and combining the two yields $S\vdash \varphi$. Hence, $g$ is faithful.
\end{proof}
\begin{lemma}[Upwards Preservation of Rigidity for (Co-)Retractions]\label{upwards_preservation_retractions}
     If $S$ is rigid with respect to (co-)retractions and $S\subseteq T$ is an extension of $S$ in the same finite language, then $T$ is rigid with respect to (co-)retractions if $T$ is finitely axiomatized relative to $S$. 
\end{lemma}
\begin{proof}
    Let $T=S+\psi$ and let $f:T\to T$ be a retraction. Then there is an interpretation $g:T\to T$ such that $f\circ g$ is homotopic to the identity over $T$ via a homotopy $\eta$. Because the language of $T$ is finite, this can be expressed in a sentence $\Iso(\eta,f\circ g,\id_T)$ proved by $T$. Define disjunctive interpretations $f':S\to S$ and $g':S\to S$ such that $S\vdash f'(\varphi)\leftrightarrow((\psi\land f(\varphi))\lor(\neg\psi\land \iota_{d_f}(\varphi)))$ and $S\vdash g'(\varphi)\leftrightarrow((\psi\land g(\varphi))\lor(\neg\psi\land \iota_{d_g}(\varphi)))$ for $\mathcal{L}$-formulae $\varphi$. Applying $f'$ to the second equivalence and using that $S\vdash f'(\psi)\leftrightarrow \psi$, one can show that $S\vdash (f'\circ g')(\varphi)\leftrightarrow((\psi\land (f\circ g)(\varphi))\lor(\neg\psi\land \iota_{d_fd_g}(\varphi)))$ for all $\mathcal{L}$-formulae $\varphi$, that is, that $f'\circ g'$ is the disjunctive interpretation for $f\circ g$ and $\iota_{d_fd_g}$ over $\psi$. Now, one can define a disjunctive homotopy $\eta'(\bar x_0,y):=(\psi\land \eta(\bar x_0,y))\lor (\neg \psi\land  \theta_{d_fd_g}(\bar x_0,y))$, where $\theta_{d_fd_g}(\bar x_0,y):=\Delta_{d_fd_g}(\bar x_0)\land x_0=y$. It is easy to check by cases that $S\vdash \Iso(\eta',f'\circ g', \id_S)$. Therefore, $f'$ is a retraction, and the remainder of the proof is identical to the proof of Lemma~\ref{upwards_preservation_rigidity}. The case for co-retractions is analogous.
\end{proof}
\begin{lemma}[Backwards Preservation of Rigidity under Retracts]\label{rigid_under_retract}
If $S$ is $C,E$-rigid, and $T$ is a $C,E$-retract of $S$, that is, if there is $f:T\to S$ and $g:S\to T$ with $f,g\in C$ and $g\circ f\mathrel{E} \id_T$, then $T$ is $C,E$-rigid.
\end{lemma}
\begin{proof}
    Let $h:T\to T$ be an interpretation in $C$. Then $f\circ h \circ g:S\to S$ is in $C$, whereby $f\circ h\circ g \mathrel{E} \id_S$. But then $h\mathrel{E}(g\circ f)\circ h\circ (g\circ f)  \mathrel{E} g\circ (f\circ h \circ g)\circ f \mathrel{E} (g\circ f)\mathrel{E}\id_T$. 
\end{proof}

\noindent
From this lemma, we immediately get that rigidity is preserved by bi-interpretations. With the aid of the following definition, which includes sentential equivalence, synonymy, and bi-interpretability as special cases, we can generalize this observation: 
\begin{definition}[$C,E$-Equivalence]
    Two theories $S$ and $T$ are $C,E$-equivalent iff there are interpretations $f:S\to T$ and $g:T\to S$ such that $f,g\in C$ and $g\circ f\mathrel{E}\id_S$ and $f\circ g\mathrel{E}\id_T$.
\end{definition}
\noindent
It is then immediately clear that $C,E$-rigidity is preserved by $C,E$-equivalence. Finally, let us note that $C,E$-rigid theories obey a Cantor--Schröder--Bernstein theorem with respect to $C$ and $E$. More precisely:
\begin{theorem}[CSB for Rigid Theories]
    Let $S$ and $T$ be $C,E$-rigid. Assume that there are interpretations $f:S\to T$ and $g:T \to S$ such that $f\in C$ and $g\in C$. Then $S$ and $T$ are $C,E$-equivalent.
\end{theorem}
\begin{proof}
By assumption $f\in C$ and $g\in C$, and $C$ is closed under composition, so the
self-interpretations $g\circ f:S\to S$ and $f\circ g:T\to T$ both lie in $C$. Since $S$ is
$C,E$-rigid and $g\circ f\in C$, we have $g\circ f\mathrel{E}\id_S$; likewise, $f\circ g\mathrel{E}\id_T$
because $T$ is $C,E$-rigid. Therefore, $S$ and $T$ are $C,E$-equivalent via $f$ and $g$.
\end{proof}

\section{Arithmetic}
In this section, we study various arithmetical theories with respect to their rigidity properties.
\subsection{Presburger Arithmetic}
That Presburger Arithmetic, the theory of addition on the natural numbers, is rigid was originally conjectured by Visser. A particular consequence of this conjecture is that $\PrA$ is not interpretable in any of its finitely axiomatized subtheories, a property it shares with $\PA$. The conjecture was first investigated by \cite{zoethout2015interpretations} and then settled positively by \cite{Pakhomov2020}. Using the theory of order types of non-standard models of $\PrA$ and a notion of rank for linear orders, they show that any model of $\PrA$ interpreted in the standard model must be isomorphic to the standard model. The proof that this isomorphism is furthermore definable requires a theory of definable families of sets in $\PrA$ and the eliminability of counting quantifiers over $\PrA$ established by \cite{apelt1966} and \cite{schweikardt2005arithmetic}. The theorem remains true even if one allows identity-modification and parameters.
\begin{theorem}[Rigidity of Presburger Arithmetic, \citet{Pakhomov2020}]\label{rigidity_pra}
The theory $\PrA$ is rigid. 
\end{theorem}

\subsection{Successor Arithmetic}

\begin{definition}[Successor Arithmetic]
$\SucA$ is the first-order theory in $\mathcal{L}=\{0,S\}$ with the following axioms:
\begin{itemize}
    \item[$(S1)$] $\forall x (Sx \neq 0)$
    \item[$(S2)$] $\forall x \forall y (Sx = Sy \to x = y)$
    \item[$(S3)$] $\forall x (x \neq 0 \to \exists y (x = Sy))$
    \item[$(S4^m)$] $ \forall x (S^m x \neq x)$ for each $m \geq 1$
\end{itemize}
Let $\NN^S = \langle\N, 0, S\rangle$ be the standard model of $\SucA$. 
\end{definition}

\medskip
\noindent
By \cite[Theorem 2.5.8]{hinman2005}, $\SucA = \Th(\NN^S)$; by Corollary 2.5.12, ibid, $\SucA$ is $\kappa$-categorical for $\kappa \geq \aleph_1$. The countable models can be characterized as follows:

\begin{fact}[Classification of Models]\label{SucA-models}
Every countable model of\/ $\SucA$ is isomorphic to
$\N + \lambda\times\Z$ for a cardinal~$\lambda \le \aleph_0$, where
$\N + \lambda\times\Z$ denotes the successor structure with
underlying set $\N\cup(\lambda\times\Z)$, zero element
$0\in\N$, and successor acting as $n\mapsto n+1$ on the
$\N$-part and as $(\alpha,k)\mapsto(\alpha,k+1)$ on each
$\Z$-chain.
\end{fact}
  \noindent
We recall some facts about definability in $\NN^S$. Since $\SucA$ admits quantifier elimination \citep[Ex. 2.4.53]{hinman2005}, the definable sets are precisely those cut out by quantifier-free formulae. 
 
\begin{fact}[Definable Sets]\label{SucA-sets}
The definable subsets of\/ $\N$ in $\langle\N,0,S\rangle$ are exactly the finite and cofinite sets.
\end{fact}

\begin{fact}[Definable Functions]\label{SucA-funct}
Let $D\subseteq\N$ be cofinite and definable, and let $g:  D\to\N$ be definable in $\NN^S$. Then there is $M\in\N$ such that either $g$ is constant on $\{x\in D:x>M\}$, or there is $c\in\Z$ with $g(x)=x+c$ for all $x\in D$ with $x>M$.
\end{fact}

\begin{theorem}[Rigidity of Successor Arithmetic]\label{SucA-a-rigid}
The theory $\SucA$ is rigid with respect to one-dimensional interpretations, up to homotopy.
\end{theorem}

\begin{proof}
Let $f$ be a one-dimensional self-interpretation of $\SucA$ and let $\NN^S_f=\langle\N_f,0_f,S_f\rangle$ be the interpreted model in $\NN^S$. Since $\SucA=\Th(\NN^S)$ is complete and $\mathcal{L}$ is finite, by Lemma~\ref{Iso} it
suffices to exhibit a definable isomorphism $h: \NN^S\to\NN^S_f$.

\medskip\noindent
The definable set $\N_f$ is infinite, hence cofinite by Fact~\ref{SucA-sets}. The definable function $S_f$ is injective on $\N_f$ by $(S2)$, hence not eventually
constant; so by Fact~\ref{SucA-funct} there are $m_0$ and $c\in\Z$ with $S_f(x)=x+c$ for all $x\in\N_f$ with $x>m_0$. Enlarging $m_0$, we may assume $m_0\geq\max\bigl((\N\setminus\N_f)\cup\{0_f\}\bigr)$. Put
$C:=\{x:x>m_0\}\subseteq\N_f$ and $F:=\N_f\setminus C$, a finite set containing~$0_f$.

\medskip\noindent
By $(S1)$--$(S3)$, $S_f$ is a bijection of $\N_f$ onto $\N_f\setminus\{0_f\}$, so
\[
  S_f(F)\cup S_f(C)=(F\setminus\{0_f\})\cup C,
  \qquad S_f(C)=\{m_0{+}1{+}c,\,m_0{+}2{+}c,\ldots\}.
\]
If $c\geq 0$, then $S_f(C)$ misses exactly $c$ elements of $C$; if $c<0$, then $S_f(C)$ exceeds $C$ by $|c|$ elements, which lie in $\N_f$ below $C$, i.e., in
$F$, and differ from $0_f$ by $(S1)$. In either case, comparing cardinalities along the decomposition yields $|S_f(F)|=|F|-1+c$; since $S_f$ is injective and
$F$ is finite, $|S_f(F)|=|F|$, whence $c=1$. In particular, $C$ is closed under $S_f$, which acts on $C$ as $x\mapsto x+1$.

\medskip\noindent
Set $h(k):=S_f^k(0_f)$. The orbit of $0_f$ is injective, since a repetition would produce an $S_f$-cycle, contradicting $(S4^m)$; thus $h$ is injective,
with $h(0)=0_f$ and $h(k{+}1)=S_f(h(k))$. Being infinite, the orbit leaves the finite set $F$, say $h(n_0)\in C$; as $S_f$ acts on $C$ as $x\mapsto x+1$, we get
$h(n)=n+d$ for all $n\geq n_0$, where $d:=h(n_0)-n_0$. Hence the orbit of $0_f$ is cofinite in $\N_f$. By Fact~\ref{SucA-models}, the countable model
$\NN^S_f$ is isomorphic to $\N+\lambda\times\Z$, where the $\N$-part corresponds exactly to the orbit of $0_f$ and each $\Z$-chain is infinite and disjoint from it. Since the complement of the orbit in $\N_f$ is finite,  $\lambda=0$, so $h$ is surjective and $h: \NN^S\cong\NN^S_f$.

\medskip\noindent
On the cofinite definable set $\{n:n\geq n_0\}$ the map $h$ is the shift $n\mapsto n+d$, and at the finitely many $k<n_0$ it takes the definable values $S_f^k(0_f)$; hence its graph is a finite union of definable sets and $h$ is
definable in $\NN^S$ by a formula $\eta(x,y)$. The sentence
$\Iso(\eta,\id, f)$ is true in $\NN^S$ and therefore provable in $\SucA$, so $f$ is homotopic to $\id_{\SucA}$.
\end{proof}

\medskip\noindent
\emph{Remark.}
The theorem fails for dimensions $n \geq 2$, that is, $\SucA$ is not rigid with respect to multi-dimensional interpretations; in particular, interpreted models need not be isomorphic to $\NN^S$. This is because in $\mathbb{N}^2$ there are definable sets that
are neither finite nor cofinite, such as the line $\{(x,y) : x = 0\}$. Consider the $2$-dimensional interpretation with domain $\N_f=\{(x,y) : x=0\lor x=S0\lor x=SS0\}$, $0_f=(2,0)$, and
\[
S_f(2,y)=(2,Sy),\quad S_f(1,Sy)=(1,y),\quad S_f(1,0)=(0,0),\quad S_f(0,y)=(0,Sy).
\]
The interpreted structure consists of the orbit $(2,0)\to(2,1)\to(2,2)\to\cdots$ of $0_f$ together with the $\Z$-chain $\cdots\to(1,1)\to(1,0)\to(0,0)\to(0,1)\to\cdots$, so it is isomorphic to $\N+\Z$ and hence a model of $\SucA$ that is not isomorphic to $\NN^S$.

\subsection{Skolem Arithmetic}

\begin{definition}[Skolem Arithmetic]
$\SkA$ is the first-order theory in $\mathcal{L}=\{1,\times\}$ of $\Nx= \langle \mathbb{N}_{>0},1,\times \rangle$.
\end{definition}

\medskip\noindent
Skolem Arithmetic can be given an explicit axiomatization \citep{cegielski1981}\footnote{An English summary of Cegielski's results can be found in chapter 3.5 of \citet{smorynski2012logical}.} and can be shown to be decidable via its description as a weak direct power of Presburger Arithmetic \citep{mostowski1952direct}. In more contemporary terms, the primary consequences of this description are known under the label Feferman--Vaught theorem~\citep{feferman1959first}; see \citet{bes2008application} for a contemporary presentation. 

\medskip\noindent
Although there are self-interpretations of Skolem Arithmetic which are not provably equivalent to the identity, such as the one which restricts multiplication to the squares, we will show in this subsection that $\SkA$ is rigid up to homotopy with respect to one-dimensional interpretations. Roughly, the idea is that one can use the homogeneity of the standard model $\Nx$ of $\SkA$ to first show that the new primes in any interpreted model are uniform powers of the old primes and that they induce a finitary support structure, and to then adduce the rigidity result for $\PrA$ to get definable isomorphisms between the sets of powers of the old and the new primes. Again using the homogeneity of $\Nx$, one can then aggregate these local isomorphisms into a homotopy.

\medskip\noindent
First, let us record the following fact, which will be used in almost all stages of the proof:
\begin{fact}[Automorphisms Induced by Finitely Supported Permutations of Primes]\label{prime_hom}
    Let $\Nx$ be the standard model of $\SkA$ and let $P=\{p_1,\dots,p_n\}$ and $Q=\{q_1,\dots,q_n\}$ be two finite sets of primes. Any bijection $b:P\to Q$ extends to a permutation of the set of primes fixing every prime outside of $P\cup Q$, and any permutation $\pi:\mathbb{P}\to \mathbb{P}$ of the primes extends multiplicatively to an automorphism $\pi:\Nx\to \Nx$.   
\end{fact}

\noindent
To execute this plan, we first need a formal definition of \emph{prime} and \emph{support}:
\begin{definition}[Support of a Number]\label{support}
    Let $\M$ be a model of $\SkA$. The set of all $\M$-primes $\mathbb{P}_\M$ in $\M$ is defined by the formula $\mathrm{prime}(p):\leftrightarrow p\neq 1 \land \forall x \, (x\mid p\to x=p \lor x=1)$. The support $\supp_\M(x)$ of an element $x\in M$ is defined by the formula $\supp(p,x):\leftrightarrow\mathrm{prime}(p)\land p\mid x$.
\end{definition}
\noindent
Functions definable in $\Nx$ preserve support in the following sense:
\begin{lemma}[Support Preservation under Definable Functions]\label{support preservation}
    Let $f:\SkA\to\SkA$ be a one-dimensional relative interpretation. Let $x,y\in \N^\times_f$. Then $\supp_{\Nx}(x\times_f y)\subseteq \supp_{\Nx}(x)\cup \supp_{\Nx}( y)$. 
\end{lemma}
\begin{proof}
    First note that since $\times_f$ is definable, $x\times_fy$ is fixed by any automorphism fixing each prime in $\supp(x)\cup\supp(y)$. Towards a contradiction, assume that $q\mid x\times_f y$ for a prime $q\notin \supp_{\Nx}(x)\cup \supp_{\Nx}( y)$. Take another prime $r\notin \supp_{\Nx}(x)\cup \supp_{\Nx}( y)$ and consider the automorphism $\sigma$ which switches $q$ and $r$ and fixes all other primes. Then $\sigma (x\times_f y)=x\times_f y$ and hence $r\mid x\times_f y$. But there are infinitely many such primes $r$, whence $x\times_f y$ would have to be divisible by infinitely many primes. Contradiction. 
\end{proof}

\noindent 
Note that for such an interpretation $f$, we have $1_f=1$, because there are no other definable elements in the standard model.
\begin{lemma}[New Primes Are Powers of Old Primes]
    Let $f:\SkA\to\SkA$ be a one-dimensional relative interpretation. Let $x$ be an $\Nx_f$-prime. Then there is an integer $e\geq 1$ such that $x=p^e$ for an $\Nx$-prime $p$.
\end{lemma}
\begin{proof}
    Towards a contradiction, let $x$ be an $\Nx_f$-prime with prime decomposition $p_1^{e_1} p_2^{e_2}\ldots p_n^{e_n}$ where $p_i\neq p_j$ if $i\neq j$ and $e_i\geq 1$, and assume $n\geq 2$. Then, by Fact~\ref{prime_hom}, there is another $\Nx_f$-prime $y$ with decomposition $q_1^{e_1} q_2^{e_2}\ldots q_n^{e_n}$, where $q_i\neq p_j$ for any $1\leq i,j\leq n$ and $q_i\neq q_j$ if $i\neq j$. Denote by $\tau$ the automorphism swapping $p_i$ and $q_i$ for $1\leq i\leq n$. Now consider the product $z:=x \times_f y$. By Lemma~\ref{support preservation}, $\supp(z) \subseteq \{p_1, \ldots, p_n, q_1, \ldots, q_n\}$. Since $\times_f$ is commutative (as $\Nx_f \vDash \SkA$) and $\tau(x) = y$, $\tau(y) = x$, we have $\tau(z) = \tau(x) \times_f \tau(y) = y \times_f x = z$. Since $\tau$ fixes $z$ and swaps $p_i$ and $q_i$, the exponents of $p_i$ and $q_i$ in the prime decomposition of $z$ are equal for each $i$. Therefore the automorphism $\rho$ defined by $p_1 \mapsto q_1$, $q_1 \mapsto p_1$, fixing all other primes, also fixes $z$. This implies that $\rho(x)=q_1^{e_1}p_2^{e_2}\ldots p_n^{e_n}$ and $\rho(y)=p_1^{e_1} q_2^{e_2}\ldots q_n^{e_n}$ are two new $\Nx_f$-primes with product $z$, contradicting uniqueness of decompositions into two prime factors (a first-order property) in $\Nx_f$.
\end{proof}
\noindent
We now record some basic definability facts about the structure of models of $\SkA$, following \cite{stonestrom2022some}, which will allow us to understand models of $\PrA$ embedded in $\Nx$.
\begin{definition}[Powers of a Prime and $p$-Part]\label{def:p-part}
    For an $\M$-prime $p$ in a model $\M$ of $\SkA$, define the formula
    \[
        \mathrm{pow}_p(x): \;\leftrightarrow\; \forall q\,\bigl(\mathrm{prime}(q) \land q\mid x \to q = p\bigr),
    \]
    which expresses that $\supp_\M(x)\subseteq\{p\}$. Such an $x$ is called a power of $p$. The set $M_p := \{x \in M : \mathrm{pow}_p(x)\}$ consists of all powers of $p$ in $\M$, including $1$.

\medskip\noindent
    $\SkA$ proves that for every $\M$-prime $p$ and every $x\in M$, there exists a unique $y\in M$ satisfying
    \[
        \mathrm{pow}_p(y) \;\land\; y \mid x \;\land\; \forall z\,\bigl(\mathrm{pow}_p(z) \land z \mid x \to z \mid y\bigr),
    \]
    i.e., $y$ is the greatest power of $p$ dividing~$x$, with respect to the divisibility ordering. We call this $y$ the \emph{$p$-part} of $x$ and denote it by $x_{p}$.
\end{definition}

\begin{lemma}[Models of $\PrA$ in the Standard Model of $\SkA$]\label{Stable_Embedding_of_PrA_in_SkA}
     Let $\M$ be a model of $\SkA$ and let $p$ be a prime in $\M$. Then $\M_p:=\langle M_p, 1, \times_{\restriction M_p}\rangle$ is a model of $\PrA$, with $p$ playing the additive role of $1$. Moreover, in the standard model $\Nx$, a subset $Q$ of $(\N^\times_p)^k$ for $k\geq 1$ is definable in $\Nx$ (with parameter $p$) if and only if $Q$ is definable in the language of $\PrA$ in the model $\Nx_p$ (without parameters). 
\end{lemma}
\begin{proof}

We use two facts from \cite{stonestrom2022some}, in the notation of Definition~\ref{def:p-part}. For a $\{1,\times\}$-formula $\varphi$ and a prime $p$, write $\varphi^{(p)}$ for the \emph{relativization of $\varphi$ to $p$}, obtained by replacing every variable by its $p$-part. Then:
\begin{itemize}
    \item[\textnormal{(R)}] \emph{\cite[Rem.~2.2]{stonestrom2022some}}: for every formula $\varphi(\bar v)$, every $\M\vDash\SkA$, every prime $p$ of $\M$, and every tuple $\bar a$ in $\M$,
    \begin{equation*}
        \M\vDash\varphi^{(p)}(\bar a)\iff   \M_p\vDash\varphi(\bar a_p),
    \end{equation*}   where $\bar a_p$ is the tuple of $p$-parts of $\bar a$;
    \item[\textnormal{(QE)}] \emph{\cite[Thm.~A.4]{stonestrom2022some}}: every formula $\vartheta(\bar v)$ is equivalent modulo $\SkA$ to a Boolean combination of formulae $\chi^{\geq m}(\bar v):=\exists^{\geq m}u\,(\mathrm{prime}(u)\wedge\chi^{(u)}(\bar v))$, where $\chi$ is a $\{1,\times\}$-formula and $m\in\N$.\footnote{Stonestrom presents a refinement of Mostowski's quantifier elimination result for $\SkA$ originally due to Je\v r\'abek, but we will not need more details.}
\end{itemize}

\medskip\noindent
By \cite[Rem.~2.1]{stonestrom2022some} (and completeness of $\SkA$), $\M_p=\langle M_p,1,\times\rangle\vDash\PrA$.

\medskip\noindent For the ``if'' direction: if a set $Q\subseteq(\N^\times_p)^k$ is defined in $\Nx_p$ by
$\psi(\bar v)$, then by (R) it is defined in $\Nx$ by
$\bigwedge_{i\leq k}\mathrm{pow}_p(v_i)\wedge\psi^{(p)}(\bar v)$, with parameter $p$.

\medskip\noindent
For the ``only if'' direction: let $Q\subseteq (\N^\times_p)^k$ be definable with parameter $p$ in $\Nx$. Then there is a formula $\vartheta(\bar v,x)$ with $|\bar v|=k$ such that $Q=\{\bar a\in\N_{>0}^k:\Nx\vDash\vartheta(\bar a,p)\}$. By (QE), $\vartheta(\bar v,x)$ is a Boolean combination of formulae of the form $\chi^{\geq m}(\bar v,x)$; since a Boolean combination of $\Nx_p$-definable sets is again $\Nx_p$-definable, it suffices to define the restriction of a single such formula to $(\N^\times_p)^k$.

\medskip\noindent
So fix $\chi$ and $m$, and let $\bar a\in(\N^\times_p)^k$. By (R),
\begin{equation*}
     \Nx\vDash\chi^{\geq m}(\bar a,p)
    \iff
    \bigl|\{q\in\mathbb{P}:\Nx_q\vDash\chi(\bar a_q,p_q)\}\bigr|\geq m.
\end{equation*}
Here we modify the proof of \cite[Lem.~2.16]{stonestrom2022some}. Clearly, $p_q=p$ if $p=q$ and $1$ if $p\neq q$. Furthermore, as $\bar a$ is a tuple of powers of $p$, we have $\bar a_p=\bar a$, while $\bar a_q=\bar 1$ for every prime $q\neq p$. Thus a prime $q\neq p$ contributes to the count precisely when $\Nx_q\vDash\chi(\bar 1, 1)$, a single sentence whose truth value is the same at every such $q$ and is independent of $\bar a$; the prime $p$ contributes precisely when $\Nx_p\vDash\chi(\bar a,p)$. Writing $N\in\{0,\infty\}$ for the number of contributing primes $q\neq p$ (so $N=\infty$ if $\chi(\bar 1,1)$ holds in $\langle\N,0,+\rangle$ and $N=0$ otherwise), the count is $N$ plus a $0/1$ contribution from $p$, so
\begin{equation*}
    \Nx\vDash\chi^{\geq m}(\bar a,p)
    \iff
    \begin{cases}
        \text{true for all }\bar a, & \text{if }N=\infty\text{ or }m=0,\\[2pt]
        \Nx_p\vDash\chi(\bar a,p), & \text{if }N=0\text{ and }m=1,\\[2pt]
        \text{false for all }\bar a, & \text{if }N=0\text{ and }m\geq 2.
    \end{cases}
\end{equation*}
In each case the restriction of $\chi^{\geq m}(\bar v,p)$ to $(\N^\times_p)^k$ is defined in $\Nx_p$ by one of the $\{1,\times\}$-formulae $\top$, $\bot$, or $\chi(\bar v,p)$, with parameter $p\in\N^\times_p$.

\medskip
\noindent
Taking the matching Boolean combination over the formulae $\chi^{\geq m}(\bar v,p)$ for $\vartheta$ yields a $\{1,\times\}$-formula $\psi(\bar v, p)$ with $Q=\{\bar a\in(\N^\times_p)^k:\Nx_p\vDash\psi(\bar a,p)\}$. Finally, $p$ is definable without parameters in $\Nx_p\cong\langle\N,0,+\rangle$ as the least non-trivial element in the definable ordering of $\PrA$. Replacing the parameters by their definitions in $\psi$ shows $Q$ is definable without parameters in the language of $\PrA$ in $\Nx_p$.
\end{proof}

\begin{lemma}[Inclusion of Sets of Powers]\label{inclusion_of_fibers}
    Let $f:\SkA\to\SkA$ be a one-dimensional relative interpretation. Let $p^e$ be an $\Nx_f$-prime, with $p$ an $\Nx$-prime and $e\geq 1$. Then $(\N^\times_f)_{p^e}\subseteq \N^\times_p$.
\end{lemma} 
\begin{proof}
    Let $x\in (\N^\times_f)_{p^e}$. If $x=1$, the inclusion is clear. We want to show that if $a\mid x$ and $a\neq 1$, then $p\mid a$. In $\Nx$, this is the case iff $\supp(x)=\{p\}$. Assume that $x\in (\N^\times_f)_{p^e}$ but that there is a prime $q\neq p$ such that $q\mid x$. Pick another prime $r\notin \supp(x)\cup\{ p,q\}$. Now, let $\sigma$ be the automorphism swapping $q$ and $r$. Then $\sigma(x)\in (\N^\times_f)_{p^e}$: if $a\neq 1$ and $a\mid_f\sigma(x)$, then $\sigma(a)\mid_f x$, whereby $p^e\mid_f \sigma(a)$. Since $p^e$ is fixed by $\sigma$, this implies $p^e\mid_fa$. Now, $(\Nx_f)_{p^e}$ is a model of $\PrA$, linearly ordered by $\mid_f$. So $x\mid_f\sigma(x)$ or $\sigma(x)\mid_fx$. In the first case, applying $\sigma$ yields $\sigma(x)\mid_fx$. By antisymmetry, this implies $x=\sigma(x)$. Mutatis mutandis, the second case yields the same conclusion. But this conclusion is contradicted by  the choice of $q\in \supp(x)\setminus \supp(\sigma(x))$.
\end{proof}
\begin{lemma}[New Primes Are Uniform Powers of Old Primes]\label{uniformity_of_new_primes}
    Let $f:\SkA\to\SkA$ be a one-dimensional relative interpretation. Then there is an integer $e\geq 1$ such that $x\in \N^\times_f$ is an $\Nx_f$-prime iff $x=p^e$ for a prime $p\in \N$.
\end{lemma}
\begin{proof}    
     We first show that for no prime $p\in \N$, both $p^a$ and $p^b$ are $\Nx_f$-primes if $a\neq b$. Assume that $p\in \N$ is a prime such that $p^a$ and $p^b$ are primes of $\Nx_f$. Note that $(\N^\times_f)_{p^a}\times (\N^\times_f)_{p^b}\subseteq \N^\times_p\times \N^\times_p$ by Lemma~\ref{inclusion_of_fibers}. Define a function $(\N^\times_f)_{p^a}\times (\N^\times_f)_{p^b}\to \N^\times_p$ by $(x,y)\mapsto x\times_f y$. This is a well-defined function by support preservation, and definable inside the model of Presburger Arithmetic $\Nx_p$ by Lemma~\ref{Stable_Embedding_of_PrA_in_SkA}, as it is definable in $\Nx$ with parameter $p$. Moreover, it is an injection by uniqueness of prime decompositions in $\Nx_f$. This is a contradiction: the set $(\N^\times_f)_{p^a}\times(\N^\times_f)_{p^b}$ has Presburger dimension $2$ and an injection to $\N^\times_p$ would have to grow at least quadratically, but any Presburger definable function is piecewise linear.\footnote{See \cite{Pakhomov2020} for the theory of Presburger dimensions.}

    \medskip\noindent
    Finally, using Fact~\ref{prime_hom}, if $p^e$ is an $\Nx_f$-prime for some prime $p$, then $q^e$ is an $\Nx_f$-prime for any prime $q$.  
\end{proof}
\begin{lemma}[Finite Support in Interpreted Models]\label{Finite_Support}
    Let $f:\SkA\to\SkA$ be a one-dimensional relative interpretation and let $e \geq 1$ be as above. Let $x\in \N^\times_f$. Then $p^e\mid_f x$ implies $p\mid x$.
\end{lemma}
\begin{proof}
    Let $x\in \N^\times_f$ with support $\supp(x)=\{p_1,\dots,p_n\}$ and assume that $q\in \mathbb{P}\setminus \supp(x)$ such that $q^e\mid_fx$. Then, since every bijection of the primes induces an automorphism of $\Nx$ by multiplicative extension, for any other prime $r\notin \supp(x)$, there is an automorphism $\pi:\Nx\to \Nx$ which is the identity on $\supp(x)$ and maps $q$ to $r$. Hence $r^e\mid_fx$. Therefore, the only primes not dividing $x$ in $\Nx_f$ are contained in the set $\{p_1^e,\dots,p_n^e\}$. But $\SkA$ proves each number to have infinitely many prime non-divisors, expressed schematically. Contradiction.
\end{proof}

\noindent
This lemma shows that in an interpreted model $\Nx_f$, each number has finite support in the interpreted model. This rules out one source of non-standardness. However, it might be that some of the interpreted primes behave non-standardly. More precisely, it might be that the interpreted set of powers of a prime $(\N^\times_f)_{p^e}$ gives rise to a non-standard model of $\PrA$. However, this can be ruled out using Lemma~\ref{Stable_Embedding_of_PrA_in_SkA} and Theorem~\ref{rigidity_pra}. 

\begin{lemma}[Models of Presburger Arithmetic in Interpreted Models]\label{standard_fibers}
    Let $f:\SkA\to\SkA$ be a one-dimensional relative interpretation and let $e\geq 1$ be as above. Then there is a formula $\theta(x,y,z)$ such that for each prime $p$, $\theta(x,y,p)$ defines an isomorphism $\theta_p:\Nx_p\to(\Nx_f)_{p^e}$. 
\end{lemma}
\begin{proof}
Due to Lemma~\ref{inclusion_of_fibers},
$(\N^\times_f)_{p^e}\subseteq\N^\times_p$. Moreover, by Lemma~\ref{Stable_Embedding_of_PrA_in_SkA}, the structure $(\Nx_f)_{p^e}$ is definable in $\Nx_p$, giving us an interpretation of $\PrA$ into $\PrA$. By Theorem~\ref{rigidity_pra}, there is hence a definable isomorphism $\theta_p:\Nx_p\to(\Nx_f)_{p^e}$. In $\Nx$, this isomorphism is defined by a formula $\theta(x,y,p)$ with parameter $p$, because the Presburger formula can be chosen without parameters. Invoking Fact~\ref{prime_hom} again, $\theta(x,y,q)$ defines an isomorphism between the induced models of $\PrA$ for any prime $q$, whence $\theta(x,y,z)$ satisfies the conclusion of the lemma.
\end{proof}
\noindent
We now aggregate the local isomorphisms into a global definable isomorphism and thereby establish the rigidity of Skolem Arithmetic with respect to one-dimensional interpretations.
\begin{theorem}[Rigidity of Skolem Arithmetic]\label{SkA_rigidity}
    $\SkA$ is rigid with respect to one-dimensional interpretations, up to homotopy.
\end{theorem}
\begin{proof}
    Let $f:\SkA\to\SkA$ be a one-dimensional relative interpretation. We will show that there is a homotopy $\Phi:\id\Rightarrow f$. Let $e\geq 1$ be the integer from Lemma~\ref{uniformity_of_new_primes}, so that the $\Nx_f$-primes are exactly the elements of the form $p^e$ with $p$ an $\Nx$-prime. Let $\theta(x,y,z)$ be the formula from Lemma~\ref{standard_fibers}, so that for each $\Nx$-prime $p$, $\theta(x,y,p)$ defines an isomorphism $\theta_p:\Nx_p\to(\Nx_f)_{p^e}$. We define a relation $\Phi(x,y)$ by:
    \begin{equation}\label{eq:homotopy}
       \Phi(x,y)\;\leftrightarrow\;\delta_f(y)\land \forall p\, \bigl(\mathrm{prime}(p)\to
          \theta\bigl(x_{p},\,y^f_{p^e},\,p\bigr)\bigr),
    \end{equation}
    where $x_{p}$ denotes the $p$-part of $x$ in $\Nx$, and $y^f_{p^e}$ denotes the $p^e$-part of $y$ in $\Nx_f$, i.e., the greatest power of the $\Nx_f$-prime $p^e$ dividing $y$ in the sense of $\Nx_f$.
    Functionality of $\Phi$ follows from finitude of prime decompositions in $\Nx$, Lemma~\ref{Finite_Support}, and uniqueness of prime decompositions in $\Nx_f$ (that is, from the first-order properties of prime parts, which transfer from $\Nx$ to $\Nx_f$). Totality follows by setting $y:=\prod_{i=1}^{n}{}^{\!f}\,\theta_{p_i}(x_{p_i})$ for $\supp(x)=\{p_1,\ldots,p_n\}$, where the product is taken in the sense of $\Nx_f$: by uniqueness of prime decompositions in $\Nx_f$, we then have $y^f_{p_i^e}=\theta_{p_i}(x_{p_i})$ for each $i$, while $y^f_{p^e}=1=\theta_p(1)=\theta_p(x_p)$ for every prime $p\notin\supp(x)$, so that $\Phi(x,y)$ holds.
    That $\Phi$ preserves multiplication follows from the fact that for each model $\M\vDash \SkA$, each $\M$-prime $p$ and all $x,y\in M$, $(x\times y)_p=x_p \times y_p$ and additivity of $\theta_p$ for each prime $p$. Injectivity is a consequence of injectivity of $\theta_p$ and uniqueness of prime decompositions in $\Nx$.  Finally, for surjectivity, let $y \in \N^\times_f$. By Lemmata~\ref{Finite_Support} and~\ref{standard_fibers}, the $\Nx_f$-support of $y$ is finite, and each $\PrA$-model $(\Nx_f)_{p^e}$ is isomorphic to $\Nx_p$ via~$\theta_p$. After enumerating the primes in the support of $y$, $y$ decomposes as $y = \prod_{i=1}^{k}{}^{\! f} y^f_{p_i^e}$, and each factor $y^f_{p_i^e}$ lies in the range of~$\theta_{p_i}$. Setting $x_{p_i} := \theta_{p_i}^{-1}(y^f_{p_i^e})$ for each $p_i^e$ in the $\Nx_f$-support of $y$, we obtain $x := \prod_{i=1}^k x_{p_i}$ satisfying $\Phi(x,y)$.

    \medskip\noindent Therefore, $\Phi(x,y)$ defines an isomorphism $\Nx\to\Nx_f$. By Lemma~\ref{homotopy_completeness}, $\Phi:\id\Rightarrow f$ is a homotopy between $\id$ and $f$.
\end{proof}

\medskip\noindent
    \emph{Remark.} The theorem fails for interpretations of dimension $n\geq 2$, that is, $\SkA$ is not rigid with respect to multi-dimensional interpretations up to definable isomorphism. Without investigating whether the induced models for such interpretations remain isomorphic to $\Nx$, it is clear that, in general, there will be no definable isomorphism. One can see the latter claim by considering the $2$-dimensional interpretation which defines multiplication componentwise. It is easy to see that this is indeed an interpretation, as the induced model $\Nx_f$ in $\Nx$ is a free commutative monoid on countably infinitely many generators, which also shows that it is isomorphic to the standard model. However, it is also easy to show that there is no definable isomorphism between the induced definable model and the standard model, as any such isomorphism would have to partition the primes into two definable subsets.

\subsection{Robinson Arithmetic}

Robinson Arithmetic $\QA$ is not rigid, not even with respect to faithful interpretations. This follows from a simple lemma that counts faithful self-interpretations up to sentential equivalence.

\begin{lemma}\label{fin_faith}
Let theories $S$ and $T$ be mutually faithfully interpretable. If there are only finitely many faithful self-interpretations of $S$ up to sentential equivalence, then $S$ and $T$ are sententially equivalent.
\end{lemma}
\begin{proof}
Let $f: S \to T$ and $g: T \to S$ be faithful interpretations. Then $g\circ f$ is a faithful self-interpretation of $S$. By finiteness up to sentential equivalence, there are $m> 0$ and $n\geq 0$,
such that $(g\circ f)^n$ is sententially equivalent to $(g\circ f)^{m+n}$. Therefore, for $\mathcal{L}_S$-sentences $\varphi$ 
\[
S \vdash (g\circ f)^n((g\circ f)^m(\varphi) \gdw \varphi).
\] 
Since $(g\circ f)^n$ is faithful,
\[
S \vdash (g \circ f)^m(\varphi) \gdw \varphi\]
for $\mathcal{L}_S$-sentences $\varphi$. Moreover, for $\mathcal{L}_T$-sentences $\varphi$
\[
T \vdash f \circ (g \circ f)^m \circ g (\varphi) \gdw f \circ g  (\varphi).
\]
But $f \circ g$ is a faithful self-interpretation, so
\[
T \vdash (f \circ   g)^m (\varphi) \gdw \varphi.
\]
Since $f$ is an interpretation of $S$ in $T$, $(g\circ f)^{m-1}\circ g$ is an interpretation of $T$ in $S$ and $(g\circ f)^{m-1}\circ g \circ f =(g\circ f)^m$ and $f \circ (g\circ f)^{m-1}\circ g=(f\circ g)^m$, both composites are sententially equivalent to the respective identity interpretations and $S$ and $T$ are sententially equivalent.
\end{proof}

\begin{corollary}\label{Q_faithful}
$\QA$ has infinitely many faithful self-interpretations which are pairwise not sententially equivalent, and hence pairwise non-homotopic. In particular, $\QA$ is not rigid with respect to faithful interpretations, up to homotopy.
\end{corollary}
\begin{proof}
Let $\PA^-$ be the theory of discretely ordered commutative semirings with a least element \citep[see][]{Kaye1991}. $\QA$ and $\PA^-$ are not sententially equivalent, since $\QA$ is not sententially equivalent to any sequential theory \citep[Theorem~14]{Visser2017}, while $\PA^-$ is sequential \citep{Jerabek2012}. They are mutually interpretable, and both are \emph{trustworthy} in the sense of \citet{visser2005faith}, i.e., every theory interpretable in them is faithfully interpretable in them: for $\PA^-$ this is Friedman's theorem~\citep[Corollary~5.6]{visser2005faith}, for $\QA$ it follows from \citet[Corollary~5.9]{visser2005faith}. Hence $\QA$ and $\PA^-$ are mutually faithfully interpretable. The faithful interpretations constructed in \citet[Lemma~5.4]{visser2005faith} may modify identity as a consequence of the involved Henkin construction; replacing it with the more refined Henkin interpretation of \citet[Theorem~3.15]{visser2017interpretation}, one obtains identity-preserving faithful interpretations. By Lemma~\ref{fin_faith}, $\QA$ has infinitely many faithful self-interpretations up to sentential equivalence.
\end{proof}

\noindent
In future work we will undertake a more in-depth and detailed examination of the number of self-interpretations of theories. 

\subsection{Peano Arithmetic}
From the point of view of the theory of interpretations, Peano Arithmetic has the pleasant feature that one can focus on the class of direct interpretations, at least if one is interested in only those interpretations with provably infinite relativizations. An interpretation is \emph{direct} if it is one-dimensional, unrelativized and identity-preserving.
\begin{theorem}[{\citealp[Theorem~6.13]{Visser2006}}]\label{PA_direct}
Let $T$ be an extension of $\PA$ in the language of $\PA$, and let $f:  S\to T$ be an interpretation, possibly multi-dimensional and identity-modifying, such that for each model $\M\vDash T$, the induced model $\M_f\vDash S$ is infinite. Then there is a direct interpretation $f^\circ:  S\to T$ homotopic to $f$.
\end{theorem}
\noindent
The proof only uses the least number principle in $T$, so it applies analogously to theories $T$ in extended languages $\mathcal L\supseteq\mathcal L_{\PA}$, as long as they prove the induction scheme for all $\mathcal L$-formulae. In particular, it applies to $\Th(\fA)$ for any expansion $\fA$ of $\NN$. This theorem will allow us to freely cite results from the literature about interpretability into $\PA$ without worrying too much about which features of interpretations the authors allow.

\begin{theorem}
Every consistent r.e.\ extension $T$ of $\PA$ (in the language of $\PA$) has
infinitely many pairwise non-homotopic, indeed pairwise not sententially
equivalent, self-interpretations.
\end{theorem}

\begin{proof}
Put $T_0:=T$ and $T_{n+1}:=T_n+O_{n+1}$, where $O_{n+1}$ is an Orey sentence of $T_n$, i.e., both $T_n+O_{n+1}$ and $T_n+\neg O_{n+1}$ are interpretable in $T_n$. Each $T_n$ is an r.e.\ extension of $\PA$ that is interpretable in $T$, hence consistent, and so has an Orey sentence \citep{Orey1961}. For $n\ge 1$ let
\[
  \sigma_n := O_1\wedge\dots\wedge O_{n-1}\wedge\neg O_n .
\]
Then $T+\sigma_n = T_{n-1}+\neg O_n$ is interpretable in $T_{n-1}$, and hence in
$T$. Let $f_n$ be such an interpretation of $T+\sigma_n$ in $T$; in particular,
$f_n$ is a self-interpretation of $T$ with $T\vdash f_n(\sigma_n)$.

\noindent
Let $i<j$ and assume towards a contradiction that $T\vdash f_i(\varphi)\gdw f_j(\varphi)$
for every sentence $\varphi$. From $T\vdash f_i(\sigma_i)$ we get $T\vdash f_j(\sigma_i)$;
moreover $T\vdash f_j(\sigma_j)$. On the other hand, $\sigma_i\vdash\neg O_i$ and
$\sigma_j\vdash O_i$, so $T\vdash\neg(\sigma_i\wedge\sigma_j)$, and since $f_j$ is a
self-interpretation of $T$,
\[
  T\vdash \neg\bigl(f_j(\sigma_i)\wedge f_j(\sigma_j)\bigr).
\]
Now $T$ would be inconsistent, which is a contradiction. So $f_i$ and $f_j$ are not
sententially equivalent and not homotopic either. By Theorem~\ref{PA_direct}, we may take the $f_n$ to be direct, since $T+\sigma_n$ has no finite models.
\end{proof}

\noindent
In particular, $\PA$ is not rigid. With Lemma~\ref{upwards_preservation_rigidity} it also follows that no finitely axiomatizable subtheory of $\PA$ and no subtheory that is finitely extended by $\PA$ can be rigid. We can also infer that $\PA$ is not rigid with respect to faithful interpretations. In fact, $\PA$ has an infinite number of \emph{faithful} self-interpretations. This follows from the fact that there are theories that are mutually faithfully interpretable with $\PA$, but not sententially equivalent with it, such as  $\ZF_{\mathrm{fin}}$ \citep[Theorem 5.1]{Enayat2011}. Together with Theorem~\ref{PA_direct} and Lemma~\ref{fin_faith} it follows that:

\begin{corollary}
    $\PA$ has infinitely many faithful self-interpretations, even up to sentential equivalence.
\end{corollary}
 
\medskip
\noindent
So $\PA$ is not rigid with respect to faithful interpretations. With Lemma~\ref{upwards_preservation_faithful} it also follows that no subtheory of $\PA$ that is finitely extended by $\PA$ can be rigid with respect to faithful interpretations. However, $\PA$ is rigid with respect to self-retractions and self-co-retractions. Co-retractions are a particularly robust subspecies of faithful interpretations. Visser initiated the study of (co-)retractions in the context of extensions of $\PA$. Of particular importance to us is a property that Enayat later referred to as the solidity of $\PA$: solidity says that if $\M_1,\M_2,\M_3$ are models of $\PA$ such that $\M_3$ is interpreted in $\M_2$ and $\M_2$ is interpreted in $\M_1$ and there is a definable isomorphism $i_0:\M_1\to \M_3$, then there is a definable isomorphism $i:\M_1\to \M_2$. Strictly speaking, the definition allows for parameters in the interpretations and the definitions of the relevant isomorphisms, but it is easy to see (in the case of $\PA$) that if one does not allow parameters in the setup, the resulting isomorphism will likewise be definable without parameters. Moreover, although solidity is local in the sense of being formulated at the level of interpretations of models, one can use Lemma~\ref{homotopy_compactness} to lift the resulting isomorphisms to an isomorphism between interpretations at the theory-level, if the considered languages are finite.

\begin{theorem}
$\PA$ is rigid with respect to retractions and co-retractions.
\end{theorem}
\noindent
This follows directly from \citet{Visser2006}, which establishes the solidity of $\PA$. Since solidity of a theory also transfers to consistent extensions of that theory in the same language, rigidity with respect to (co-)retractions also holds for extensions of $\PA$ (in the same language). Note that this already follows from Lemma~\ref{upwards_preservation_retractions} for finite extensions of $\PA$. 

\subsection{True Arithmetic}

In this subsection, we prove that $\Th(\NN)$ is rigid. The key insight lies in a theorem established by \citet{Feferman1958}: the standard model of arithmetic $\NN$ cannot interpret a non-standard model of True Arithmetic.\footnote{See also \cite[Theorem 25.4a]{boolos2007}. For more on this topic and its historical background, see \url{https://mathoverflow.net/questions/353898/interpreting-proper-elementarily-equivalent-end-extensions/354818}. Thanks to Ali Enayat for bringing this to our attention.} We obtain it as a corollary of a useful, more general theorem.\footnote{Many thanks to Mateusz Łełyk, Piotr Gruza, and Konstantinos Papafilippou for suggesting the more general theorem.} 

\medskip
\noindent
The following is a restatement of Feferman’s theorem~\citep{Feferman1958}, for which a simpler proof was found by Tennenbaum (according to \citet{Scott1960}). In the form presented here, it was proven by \citet[Theorem~1]{Scott1960}.

\begin{theorem}[Feferman--Scott]\label{feferman-scott}
Let $\fS = \langle \NN, (R_\xi)_{\xi<\alpha} \rangle$ be the
standard model of arithmetic $\NN = \langle \N, +, \times\rangle$ expanded by an arbitrary family of relations,
and let $\fB$ be a structure with carrier $\N$ in the same signature. If $\fB \equiv \fS$ and the interpretations in $\fB$ of $+$, $\times$ and of each $R_\xi$ are (parameter-free) definable over $\fS$, then the numeral map of $\fB$ is an isomorphism between $\fB$ and $\fS$.
\end{theorem}

\begin{theorem}\label{Th(N,A)rigid}
Let $\fA = \langle \NN, (R_\xi)_{\xi<\alpha} \rangle$ be any expansion of the standard model of arithmetic by relations (for countable $\alpha$). Then $\Th(\fA)$ is rigid.
\end{theorem}
 
\begin{proof}
Let $f$ be a self-interpretation of $\Th(\fA)$, with internal model
$\fA_f \vDash \Th(\fA)$ built inside $\fA$. Since $\fA_f$ is infinite, Theorem~\ref{PA_direct} implies that we can view it as having domain $\mathbb{N}$, with all of its operations and relations definable in $\fA$.
Moreover, since $\fA_f \vDash \Th(\fA)$ and $\Th(\fA)$ is complete, $\fA_f \equiv \fA$. Hence, Theorem~\ref{feferman-scott} applies, and its numeral map $h(n) \;=\; S_f^{\,n}(0_f)$ is an isomorphism between $\fA$ and $\fA_f$.

 \medskip
 \noindent
The graph of $h$ is definable over $\fA$ by the usual $\beta$-function coding of the finite iteration $0_f, S_f(0_f), \dots, S_f^{\,n}(0_f)$; fix a formula $\eta(x,y)$ defining it. Since $\Th(\fA)$ is complete and $\eta$ defines an isomorphism $h:\fA\to \fA_f$, $\eta:\id\Rightarrow f$ is a homotopy by Lemma~\ref{homotopy_completeness}. Therefore $\Th(\fA)$ is rigid.
\end{proof}

\begin{corollary}
$\Th(\NN)$ is rigid.
\end{corollary}

\medskip \noindent
These results also suggest a certain hierarchy with respect to the extension of complete arithmetic theories over the natural numbers. In this hierarchy, $\PrA$ occupies an initial position; unlike $\SucA$ and $\SkA$, it is rigid in all dimensions. If one extends from here to include successor and/or multiplication and considers the respective complete theories, rigidity is preserved in all dimensions and, from this point on, is no longer lost even for arbitrary relations $R$. For further investigation, it would therefore be interesting to determine whether there exist structures $\langle \N, R \rangle$ or $\langle \N, +, R \rangle$—respectively, prior to $\langle \N, + \rangle$ or between $\langle \N, + \rangle$ and $\langle \N, +, \times \rangle$—whose complete theories are rigid in all dimensions.

\medskip \noindent
In the specific example of Büchi Arithmetic $\BA_n$ for $n\geq 2$ (the theory of $\langle \N,+,V_n \rangle$, where $V_n(0)=0$ and, for $x>0$, $V_n(x)$ denotes the largest power of $n$ dividing $x$) it can be said, at least, that interpreted models are always isomorphic to the standard model; however, whether the isomorphism is definable in $\BA_n$ remains an open question \citep{ZapThesis}. 

\section{Algebra}

In this section, we consider rigid algebraic and semi-algebraic theories. Many theories in this class appear to be rigid. However, our investigations to date are still far from providing a systematic overview of this area. The selection of theories to be considered here can therefore be understood merely as a fragmentary one for this subject area. Drawing on established literature on the model theory of fields, we have initially focused on the two highly relevant field theories $\RCF$ and $\ACF$. \cite{bays2019definability}, in Section 4.2 of their paper, discuss a more model-theoretic version of what we call rigidity under the heading ``self-recollecting theories'', showing that certain theories of (valued) fields are rigid. We expand on their results mostly by eliminating parameters and establishing uniformity. In the future, the further examples regarding valued fields they mention might merit further study. One might also want to investigate perfectoid fields, whose model theory has been an exciting new area of investigation in the past few years \citep[see, e.g.,][]{jahnke2025beyond, ketelsen2026understanding}. The theory of Boolean groups $\F_2$ was chosen as an example of a non-trivially incomplete but (partially) rigid theory from a domain that model theorists have not investigated with an eye on self-interpretations.

\medskip
\noindent
\subsection{Boolean Groups}
The theory of groups is clearly not rigid, as the trivial group can be defined in any group. However, there are incomplete extensions of group theory which are rigid (with respect to one-dimensional interpretations).
\begin{definition}[Boolean Groups]
$\F_2$ is the first-order theory in $\mathcal{L}=\{0, + \}$ with the following axioms:
\begin{itemize}
\item[$(B1)$] $\forall x \forall y \forall z ((x+y)+z=x+(y+z))$
\item[$(B2)$] $\forall x \forall y (x+y=y+x)$
\item[$(B3)$] $\forall x (x+0=x)$
\item[$(B4)$] $\forall x (x + x = 0) $
\item[$(B5)$] $\exists x (x \neq 0)$
\end{itemize}
\end{definition}
\medskip\noindent
Note that any model of $\F_2$ is a vector space over the field with two elements. Therefore, given any such group, its automorphism group acts transitively on the non-zero elements, as any non-zero vector can be extended to a basis and any bijection between bases extends to a corresponding automorphism of the vector space.
\begin{theorem}
$\F_2$ is rigid with respect to one-dimensional interpretations, up to homotopy.
\end{theorem}
\begin{proof}
    Let $f$ be a one-dimensional self-interpretation of $\F_2$. In any model $\mathfrak{G}$ of $\F_2$, the interpretation $f$ specifies a domain $D \subseteq G$, a distinguished element $0_f$ and an operation $+_f:D^2\to D$ such that $\mathfrak{G}_f=\langle D,0_f,+_f\rangle$ is a model of $\F_2$.
     We distinguish two cases: If $|G|=2$, there are two possibilities for $\mathfrak{G}_f$. Let $G=\{0,a\}$. Either $0_f=0$, in which case $a_f=a$ and $+_f=+$, or $0_f=a$, in which case $a_f=0$ and the transposition $x\mapsto x+a$ is an isomorphism $\mathfrak{G}\to\mathfrak{G}_f$.
     
\medskip\noindent
If $|G|\geq 3$, we proceed as follows: By the above observation, the only definable element of $G$ is $0$, from which follows $0_f=0$. Moreover, again by the above remark, the only non-empty definable subsets of $G$ are $\{0\}$, $G\setminus \{0\}$, and $G$ itself. Non-triviality (B5) rules out $D=\{0\}$, and $0_f=0\in D$ rules out $G\setminus \{0\}$. Hence, the domain $D$ is all of $G$. 

\medskip\noindent
Now, consider two elements $a,b$ of $D$ and their sum $a+_fb$. There are three cases: In the first, $a=0_f$ or $b=0_f$. By the identity axiom~(B3) in $\mathfrak{G}_f$ and $0_f = 0$, this implies $a+_fb=a+b$. In the second, $a=b\neq 0_f$. Then, by (B4) in $\mathfrak{G}_f$, $a+_fb=0=a+b$. Finally, in the third case $a\neq b$ and $a\neq 0_f$ and $b\neq 0_f$. Now, because $+_f$ is definable, it follows that $a+_fb$ is fixed by all automorphisms fixing $a$ and $b$. From this it follows that $a+_fb \in \Span\{a,b\} =\{0,a,b,a+b\}$, the linear span of $a$ and $b$. Indeed, if $\dim(\mathfrak{G}) \leq 2$, then $G = \Span\{a, b \}$ and containment is trivial. If $\dim(\mathfrak{G}) \geq 3$, note that every element of $\Span\{ a, b \}$ is fixed by all automorphisms fixing $a$ and $b$. Conversely, if $c \notin \Span\{ a, b \}$, then $\{a, b, c\}$ is linearly independent and can be extended to a basis; replacing $c$ by $c + a$ in this basis yields a second basis, and the induced automorphism fixes $a$ and $b$ but moves $c$. Since $a +_f b$ is fixed by all automorphisms fixing $a$ and $b$, $a +_f b \in \Span\{ a, b \}$. Thus, there are four cases to consider for $a+_f b$. In the first, $a+_fb=0_f$; by (B4), $a=b$, contradicting $a\neq b$. In the second and third, namely $a+_fb=a$ and $a+_fb=b$: adding $a$ resp.\ $b$ on both sides and applying (B1) and (B4) yields $b=0_f$ resp.\ $a=0_f$, a contradiction. Therefore, $a+_fb=a+b$. Therefore $\mathfrak{G}_f=\mathfrak{G}$.

\medskip
\noindent
We can now produce a homotopy $\id\Rightarrow f$ as follows. Let $\zeta(x)$ be the formula serving as the translation of the constant $0$. Define
\begin{equation*}
\eta(x,y) \;:=\; \exists z\,\bigl(\zeta(z)\wedge y=x+z\bigr),
\end{equation*}
i.e., $\eta$ is translation by $0_f$. In every model $\mathfrak{G}$ of $\F_2$, $\eta$ defines an isomorphism $\mathfrak{G}\to\mathfrak{G}_f$: if $0_f=0$, which is always the case when $|G|\geq 3$, then $\mathfrak{G}_f=\mathfrak{G}$ and $\eta$ is the identity; in the remaining case $G=\{0,a\}$ and $0_f=a$, $\eta$ is the transposition $x\mapsto x+a$. Hence $\F_2$ proves that $\eta$ is a homotopy between $\id$ and $f$.
\end{proof}

\noindent
\emph{Remark.} The proof uses that the domain of $\mathfrak{G}_f$ is a subset of $G$, whose only definable subsets are $\emptyset$, $\{0\}$, $G\setminus\{0\}$, and $G$ when $|G|\geq 3$. For $d\geq 2$ the domain may be all of $G^d$, and $G^d$ has further definable subsets, e.g., the diagonal. For instance, the $2$-dimensional interpretation with domain $G\times G$ and componentwise addition yields $\mathfrak{G}\times\mathfrak{G}$. For finite $G$ this is not even isomorphic to $\mathfrak{G}$. For infinite $G$ it is isomorphic to $\mathfrak{G}$, but not definably so, even with parameters $\bar p$: as in the proof above, a definable map $h:  G\times G\to G$ sends $(x,y)$ into $V:=\Span\{x,y,\bar p\}$, so it maps the finite set $V\times V$ into $V$ and cannot be injective.

\medskip\noindent
\subsection{Real Closed Fields}
Moving on to a more substantial example, one can use established results on the interpretability of fields in $\RCF$ to prove the following theorem:
\begin{theorem}\label{thm:rcf}
$\RCF$ is rigid.
\end{theorem}

\begin{proof}
Let $f$ be a self-interpretation of $\RCF$. In the standard model $\RR = \langle \R, +, \times , 0, 1, < \rangle$ of $\RCF$, the interpretation $f$ specifies a domain $D \subseteq \R^n$, operations $+_f, \times_f:D^2 \to D$, distinguished elements $0_f, 1_f \in D$, an order $<_f$ on $D$, making
$\RR_f = \langle D, +_f, \times_f, 0_f, 1_f, <_f \rangle$ a model of $\RCF$.

\medskip\noindent
The real field $\RR$ is trivially an o-minimal expansion of itself, so~\cite[Theorem 1.1]{Pillay1996} applies: every infinite definable ring without zero divisors is definably isomorphic (possibly with parameters) to $\RR$, to $\RR(\sqrt{-1})$, or to the quaternions over $\RR$. Since $\RR_f$ is a commutative field, it is definably isomorphic to $\RR$ or to $\RR(\sqrt{-1})$. And since $\RR_f \vDash \RCF$, $-1$ is not a sum of squares in $\RR_f$, so $\RR_f$
is not algebraically closed. This rules out $\RR_f \cong \RR(\sqrt{-1})$. Hence, there is a definable isomorphism of rings $h:  \RR \to \RR_f$, defined by a formula $\psi(x,\bar y,\bar c)$, possibly with parameters $\bar c$ in $\RR$. Since the ordering is definable from the field structure, it is also an isomorphism of ordered fields. 

\medskip\noindent
We now show that the isomorphism is, in fact, definable without parameters. Suppose $h$ and $h'$ are two such definable isomorphisms (for possibly different parameters). Then $h^{-1} \circ h' : \RR \to \RR$ is an automorphism of $\RR$. Since $\RR$ is rigid, $h = h'$. Therefore, $\eta(x,\bar y):=\exists \bar z (\Iso(\psi,\id,f)(\bar z)\land\psi(x,\bar y,\bar z))$ defines the unique isomorphism $h:\RR\to \RR_f$ without parameters. By Lemma~\ref{homotopy_completeness}, $\eta:\id\Rightarrow f$ is the required homotopy.
\end{proof}
\noindent
\subsection{Algebraically Closed Fields}
Turning to algebraically closed fields, let us first state a crucial theorem due to Poizat, achieved through an analysis of algebraic groups definable in algebraically closed fields.
\begin{theorem}[\cite{poizat1988mm}]\label{poizat_theorem}
    If $L$ is an infinite field definable (possibly with parameters) in an algebraically closed field $K$, there is a $K$-definable isomorphism $K\cong L$ (possibly with parameters).
\end{theorem}
\noindent
With this result at hand, one can easily show that theories of algebraically closed fields are rigid, the only subtleties being uniformity and parameter-freeness of the required definable isomorphisms.\footnote{If one wanted to extend the result to identity-modifying interpretations, one could make use of elimination of imaginaries for algebraically closed fields \citep{poizat1983theorie}. We also direct the reader to \citet{poizat2026parametres}, where a deeper study of definitions with parameters in algebraically closed fields is conducted.} In the following, $\ACF_p$ denotes the theory of algebraically closed fields of characteristic $p$, where $p$ can be $0$ or any prime number. $\ACF$ denotes the theory of algebraically closed fields of arbitrary characteristic. 
\begin{theorem}
    For $p=0$ and every prime $p$, $\ACF_p$ is rigid. Moreover, the theory $\ACF$ is rigid. 
\end{theorem}
\noindent
\begin{proof}    
Let $p$ be $0$ or a prime, let $f:\ACF_p\to\ACF_p$ be an $n$-dimensional interpretation, and let $\mathfrak{K}\vDash \ACF_p$ be an algebraically closed field. By completeness of $\ACF_p$, it suffices to show that $\mathfrak{K}_f$ is definably isomorphic to $\mathfrak{K}$. Since $\mathfrak{K}_f$ is infinite, Theorem~\ref{poizat_theorem} implies the existence of an isomorphism $h:\mathfrak{K}\to  \mathfrak{K}_f$, definable with parameters by a formula $\psi(x,\bar y,\bar z)$. We now show that $h$ is definable without parameters by distinguishing cases. In case $p=0$, we first show that there can only be one isomorphism $h:\mathfrak{K}\to\mathfrak{K}_f$ definable with parameters: Letting $h_1,h_2$ denote two such isomorphisms, the fixed field $H=\{x\in K\mid (h_1^{-1}\circ h_2)(x)=x\}$ of the automorphism $h_1^{-1}\circ h_2$ is definable with parameters and hence finite or cofinite \citep[Cor.~3.4]{marker2000}. Since $H$ contains $\mathbb{Q}$, it is infinite and hence cofinite. If its complement were non-empty, it would have to be infinite, by closure under translation by elements of $H$. Therefore, $H=K$ and $h_1^{-1}\circ h_2=\id_{\mathfrak{K}}$. Then $\eta(x,\bar y):=\exists \bar z (\Iso(\psi,\id,f)(\bar z)\land\psi(x,\bar y,\bar z))$ defines a homotopy $\eta:\id\Rightarrow f$.

\medskip\noindent
Now, in case $p>0$ is prime, we no longer have unique definable isomorphisms, but the strategy is otherwise similar. Using completeness, assume that $\mathfrak{K}=\overline{\mathbb{F}_p}$. By Theorem~\ref{poizat_theorem}, we have an isomorphism $h:\mathfrak{K}\to \mathfrak{K}_f$, defined by $\psi(x,\bar y, b)$ with a single parameter $b$. A single parameter suffices because, given multiple parameters, one can replace them by polynomials in the single element $b$ generating the field extension of $\mathbb{F}_p$ generated by the parameters. Now, let $m_b(z)$ be the minimal polynomial of $b$ over $\mathbb{F}_p$. We claim that $\eta(x,\bar y)=\exists z ((m_b(z)=0)\land \psi(x,\bar y, z))$ defines a homotopy from $\id$ to $f$. Letting $F:\mathfrak{K}\to \mathfrak{K}$ denote the Frobenius $x\mapsto x^p$ and $r$ the degree of the field extension $\mathbb{F}_p(b):\mathbb{F}_p$, we see that the $r$ distinct roots of $m_b(z)$ are exactly $b,F(b),\ldots, F^{r-1}(b)$. We now need to show that $\psi(x,\bar y, F^i(b))$ defines the same isomorphism for $0\leq i\leq r-1$. Recalling that $f$ induces a functor between categories of models, let $F_f:\mathfrak{K}_f\to \mathfrak{K}_f$ denote the automorphism induced by componentwise application of $F$. This automorphism is compatible with $h$: Since $F^r$ fixes $b$, it also fixes the set defined by $\psi(x,\bar y,b)$. Therefore, $F_f^r\circ h=h\circ F^r$, which implies $(h^{-1}\circ F_f\circ h)^r=F^r$. Since $\mathrm{Aut}(\overline{\mathbb{F}_p})\cong\widehat{\mathbb{Z}}$ is abelian and torsion-free, this implies $h^{-1}\circ F_f\circ h=F$ \citep[Example 7.17]{milneFT}. Fixing some $0< i\leq r-1$, we see that the set defined by $\psi(x,\bar y,F^i(b))$ is a function given by $F_f^i\circ h\circ F^{-i}$, which, by compatibility, is $h$. 

\medskip\noindent
Rigidity of $\ACF$ now follows: By Theorem~\ref{poizat_theorem}, any self-interpretation $f:\ACF\to \ACF$ also induces a self-interpretation $f:\ACF_p\to\ACF_p$ for any $p=0$ or prime $p$. Hence, by the above, for any model $\mathfrak{K}\vDash \ACF$, there is a parameter-free definable isomorphism $\mathfrak{K}\to \mathfrak{K}_f$. Finitude of the language of rings allows one to conclude via Lemma~\ref{homotopy_compactness}. 
\end{proof}

\section{Set Theory}

Set theory partially mirrors arithmetic. It is clear that no consistent r.e.\ extension $T$ of $\ZF$ can be rigid. In fact, it must have an infinite number of non-trivial self-interpretations. This follows from the same argument as for $\PA$.\footnote{Another strategy for proving the failure of rigidity was suggested by Ali Enayat.  An \emph{inductive sequential theory} is a sequential theory that can prove the full scheme of induction over the natural numbers. Using the arithmetical completeness theorem, Tarski's undefinability of truth, and the fact that every inductive sequential theory $T$ proves the formal consistency of every finite fragment of itself, one can infer that every r.e.\ consistent inductive sequential theory fails to be rigid. Such theories include $\PA, \mathrm{ACA}, \mathrm{Z_2}, \Ze, \ZF$, and $\mathrm{KM}$.\label{Enayat_footnote}} However, similar to the case of $\PA$, it follows from the solidity of $\ZF$ \citep[Theorem 2.6]{Enayat2016} that extensions of $\ZF$ (in the same language) are rigid with respect to (co-)retractions.

\begin{theorem}
Let $T$ be an extension of $\ZF$ in the language of $\ZF$. Then $T$ is rigid with respect to (co-)retractions.
\end{theorem} 

\noindent
A key difference from arithmetic becomes apparent, however, when one considers—as in the case of $\PA$—whether there exists a complete extension of $\ZF$ or $\ZFC$ that is rigid even with respect to relative interpretations. The same applies to the question of whether the complete theory of rigid models of $\ZF$ or $\ZFC$ is itself rigid. For $\ZFC$, the following negative result can be observed; for $\ZF$ the question remains open (see Open Question below).

\begin{theorem}\label{non-rigid_ZFC}
If $T$ is a consistent extension of $\ZFC$ (in the same language), then $T$ is non-rigid.
\end{theorem}
\noindent
Let us state the idea behind the theorem. It is a well-known phenomenon that the ultrapower of a well-founded model of set theory induced by a non-principal ultrafilter on $\omega$ is always non-well-founded, so they are not isomorphic. Meanwhile, by {\L}o\'s's theorem, the ultrapower is elementarily equivalent to the original. If we internalize this construction, we can define a syntactic interpretation based on it that will yield a non-trivial self-interpretation.

\medskip\noindent
A first natural idea would be to carry out such a construction in $\ZF + V=L$, since a non-principal ultrafilter on $\omega$ is generally not definable without parameters in $\ZF$. Its existence is not provable in $\ZF$ alone; it follows from the Axiom of Choice (indeed already from the Boolean prime ideal theorem). Moreover, the constructible universe has a parameter-free definable global well-ordering that allows one to choose a canonical ultrafilter. The present, more general construction goes back to Ali Enayat, who pointed out that this should already be feasible in $\ZFC$ if, instead of using a single ultrapower of the universe, we use a direct limit of a definable collection of ultrapowers (an appropriate iterated ultrapower), as in the construction by \citet{KanoveiShelah2004}. For every consistent extension $T$ of $\ZFC$ this construction yields a self-interpretation $f$ such that for every model $\M$ of $T$ there is no $\M$-definable isomorphism between $\M$ and $\M_f$. From the point of view of $\M$, the finite ordinals of $\M_f$ properly extend the image of $\omega^{\M}$, i.e., $\M$ can ``see'' that $\M_f$ has more finite ordinals than itself.\footnote{We are very grateful to Ali Enayat for his generous explanations on this topic.}

\subsubsection*{The Kanovei--Shelah Ultrafilter}

We will repeat the ultrafilter construction from \citet[§1]{KanoveiShelah2004} for our context, since the method does not appear to be widely known. 
Everything stated can be carried out provably in $T \supseteq \ZFC$.

\medskip\noindent
Ultrafilters on $\omega$ need not admit a definable linear order, but their
\emph{enumerations} always do, because the domain $2^{\aleph_0}$ of an
enumeration is well-ordered. Let
\[
  A \;:=\; \{\, a : 2^{\aleph_0} \to \mathcal{P}(\omega) \;\mid\;
            \mathrm{ran}(a) \text{ is an ultrafilter on } \omega \,\},
  \qquad U_a := \mathrm{ran}(a),
\]
and order $A$ lexicographically: $a <_A b$ iff $a(\xi) <_{\mathrm{lex}} b(\xi)$ for the least $\xi < 2^{\aleph_0}$ with $a(\xi) \neq b(\xi)$, where for distinct $u, v \subseteq \omega$ we set $u <_{\mathrm{lex}} v$ iff
$\min(u \mathbin{\triangle} v) \in v$. Then $A$ is a set and $<_A$ is a
definable linear order on $A$. Moreover, every ultrafilter $U$ on $\omega$ satisfies $|U| \leq 2^{\aleph_0}$, so every ultrafilter occurs as $U_a$ for some $a \in A$.

 \medskip\noindent
For a finite $u \subseteq A$ and $X \subseteq \omega^A$, say $X$ is
\emph{concentrated on} $u$ in case
\[
\forall x, y \in \omega^A\,
 \bigl(x \restriction u = y \restriction u \rightarrow
  (x \in X \gdw y \in X)\bigr),
\]
and let $\mathcal{X}$ be the collection of all $X \subseteq \omega^A$
concentrated on some finite $u \subseteq A$ (the \emph{finite-support sets}).
$\mathcal{X}$ is an algebra of subsets of $\omega^A$. If $X$ is concentrated
on finite $u$ and on finite $v$, it is concentrated on $u \cap v$: given
$x \restriction (u \cap v) = y \restriction (u \cap v)$, any
$z \in \omega^A$ with $z \restriction u = x \restriction u$ and
$z \restriction v = y \restriction v$ witnesses
$x \in X \gdw z \in X \gdw y \in X$. Hence a support of minimal cardinality
is contained in every support, and each $X \in \mathcal{X}$ has a least
support $\|X\|$. For $X$ concentrated on $u$ put
$X \mathord{\downarrow} u := \{x \restriction u : x \in X\}$, so that
$x \in X \gdw x \restriction u \in X \mathord{\downarrow} u$; for
$Y \subseteq \omega^u$ and $u \subseteq v$ put
$Y \mathord{\uparrow} v := \{t \in \omega^v : t \restriction u \in Y\}$.

 \medskip\noindent
By recursion on $|u|$, define, for finite $u \subseteq A$, an iterated
ultrafilter $W_u$ on $\omega^u$: $W_\emptyset := \{\omega^\emptyset\}$, and
for $u \neq \emptyset$ with $<_A$-greatest element $b$ and
$u' := u \setminus \{b\}$,
\[
  X \in W_u \;:\gdw\;
  \{\, k \in \omega \;:\; X_{b=k} \in W_{u'} \,\} \in U_b,
  \qquad\text{where } X_{b=k} :=
  \{\, s \in \omega^{u'} : s \cup \{(b,k)\} \in X \,\}.
\]
This is the finite iteration (Fubini product) of the ultrafilters
$U_{a}$, $a \in u$, taken in $<_A$-order; cf.\
\citet[\S 6.5]{ChangKeisler1990}. Finally, let
\[
  W \;:=\; \{\, X \in \mathcal{X} \;:\;
            X \mathord{\downarrow} \|X\| \in W_{\|X\|} \,\}.
\]
 
\begin{lemma}\label{KS-ultrafilter}
$\ZFC$ proves:
\begin{enumerate}
\item[(i)] For every finite $u \subseteq A$, $W_u$ is an ultrafilter on
  $\omega^u$.
\item[(ii)] For finite $u \subseteq v \subseteq A$ and $X \subseteq \omega^u$:
  \ $X \in W_u \gdw X \mathord{\uparrow} v \in W_v$.
\item[(iii)] $W$ is an ultrafilter in the algebra $\mathcal{X}$, and for every $X\in \mathcal{X}$ and
  finite $u \supseteq \|X\|$: \ $X \in W \gdw X \mathord{\downarrow} u \in W_u$.
\end{enumerate}
\end{lemma}
 
\begin{proof}This is the finite-support ultrafilter of \citet[Proposition 2]{KanoveiShelah2004}, stated there without proof; we indicate the verifications. Part~(i) is the iterated (Fubini) product of the $U_a$, $a \in u$, in $<_A$-order, an ultrafilter by the standard induction on $|u|$  \citep[cf.][\S 6.5]{ChangKeisler1990}; parts~(ii) and~(iii) are routine inductions on $|u|$ and on a common refinement of supports, using only that each $U_b$ is an ultrafilter together with the coherence identity $(X \mathord{\downarrow} u) \mathord{\uparrow} (u \cup u')  = X \mathord{\downarrow} (u \cup u')$. The essential point in~(ii) is that the order of iteration is fixed by the single global order $<_A$ and is therefore stable under restriction. In the induction on $|v|$, one distinguishes whether the $<_A$-greatest element of $v$ belongs to $u$ or not.
\end{proof}

\subsubsection*{The Interpretation}

Let $\mathbb{F}$ be the class of all functions $g$ with
$\mathrm{dom}(g) = \omega^A$ that are concentrated on some finite
$u \subseteq A$ (i.e., $x \restriction u = y \restriction u \rightarrow g(x) = g(y)$); as before
each $g \in \mathbb{F}$ has a least support $\|g\|$. Each $g\in\mathbb{F}$ is a set. Thus, $\mathbb{F}$ is a definable proper class whose members are sets, so the quantifiers
$\exists g \in \mathbb{F} \quad\text{and}\quad \forall g \in \mathbb{F}$ appearing in the following (in $\delta_f$, in $f(x \in y)$, and in Lemma~\ref{los}) are legitimate first-order quantifiers over sets. For $g_1, \ldots, g_k \in \mathbb{F}$ and a formula $\varphi(x_1, \ldots,
x_k)$ write
\[
  \llbracket \varphi(g_1, \ldots, g_k) \rrbracket \;:=\;
  \{\, x \in \omega^A : \varphi(g_1(x), \ldots, g_k(x)) \,\},
\]
a set concentrated on $\|g_1\| \cup \dots \cup \|g_k\|$, hence in $\mathcal{X}$. In particular set
\[
g \sim_W h :\gdw \llbracket g = h \rrbracket \in W \text{ and } g \in_W h :\gdw \llbracket g \in h \rrbracket \in W.
\]

\medskip\noindent
In $\ZFC$ there may be no definable choice of representatives of $\sim_W$-equivalence classes, so we take Scott classes instead; for $g \in \mathbb{F}$ let
\[
  \hat g \;:=\; \{\, h \in \mathbb{F} \;:\; h \sim_W g \;\wedge\;
   \forall h' \in \mathbb{F}\,
   (h' \sim_W g \rightarrow \mathrm{rank}(h) \leq \mathrm{rank}(h')) \,\},
\]
a non-empty set, with $\hat g = \hat h \gdw g \sim_W h$. Set the domain $D := \{\hat g : g \in \mathbb{F}\}$. The syntactic interpretation $f$ can then be given by
\begin{align*}
\delta_f(x) \;&:=\; \exists g \in \mathbb{F}\; (x = \hat g),\\
f(x \in y) \;&:=\; \exists g, h \in \mathbb{F}\,
  \bigl( g \in x \land h \in y \land g \in_W h \bigr).
\end{align*}
The interpretation $f$ preserves identity, and note that $f(x \in y)$ does not depend on the choice of $g \in x$, $h \in y$ for $x, y$ in $D$, since $g' \sim_W g$, $h' \sim_W h$ and closure of $W$ under intersections give
$\llbracket g \in h \rrbracket \in W \gdw \llbracket g' \in h' \rrbracket \in W$.

\begin{lemma}[Internal \L o\'s]\label{los}
For each formula $\varphi(x_1, \ldots, x_k)$,
\[
  T \vdash \forall g_1, \ldots, g_k \in \mathbb{F}\,
  \bigl( f(\varphi)(\hat g_1, \ldots, \hat g_k)
   \;\gdw\;
   \llbracket \varphi(g_1, \ldots, g_k) \rrbracket \in W \bigr).
\]
\end{lemma}

\begin{proof}
The proof follows the same pattern as in \citet[Lemma 5]{KanoveiShelah2004} and proceeds by induction on $\varphi$. 

\medskip\noindent
    The atomic equality case uses $\hat g_1 = \hat g_2 \gdw g_1 \sim_W g_2 \gdw
\llbracket g_1 = g_2 \rrbracket \in W$ (by the Scott construction). The atomic membership case holds by definition of $f(x \in y)$ and the representative-independence noted above. The cases for the connectives use that
$\llbracket \neg \psi \rrbracket = \omega^A \setminus
 \llbracket \psi \rrbracket$ and
$\llbracket \psi \wedge \chi \rrbracket =
 \llbracket \psi \rrbracket \cap \llbracket \chi \rrbracket$. Since all the sets
involved lie in $\mathcal{X}$, the cases follow from $W$ being an ultrafilter in
the algebra $\mathcal{X}$ (Lemma~\ref{KS-ultrafilter}(iii)). For the existential case, let $\varphi = \exists y\, \psi(y, x_1, \dots, x_k)$, so $f(\varphi)(\hat g_1, \dots, \hat g_k)$ is $\exists y\,(\delta_f(y) \land f(\psi)(y, \hat g_1, \dots, \hat g_k))$.

\medskip\noindent
($\Rightarrow$:) Any witness is of the form $\hat h$, $h \in \mathbb{F}$; by
the induction hypothesis we have that $\llbracket \psi(h, g_1, \dots, g_k) \rrbracket \in W$ and $\llbracket \exists y\, \psi(y, g_1, \dots, g_k) \rrbracket \supseteq \llbracket \psi(h, g_1, \dots, g_k) \rrbracket$ lies in $\mathcal{X}$, hence in $W$
by upward closure within $\mathcal{X}$.

\medskip\noindent
($\Leftarrow$:) Suppose
$S := \llbracket \exists y\, \psi(y, g_1, \dots, g_k) \rrbracket \in W$ and let
$u := \|g_1\| \cup \dots \cup \|g_k\|$, on which $S$ is concentrated. By
Collection applied over the set $S$, there is a set $B$ such that
$\forall x \in S\, \exists y \in B\, \psi(y, g_1(x), \dots, g_k(x))$. Fix a well-order $\lhd$ of $B$ and define $h : \omega^A \to V$ by
\[
  h(x) :=\begin{cases}
    \text{the } \lhd\text{-least } y \in B
    \text{ with } \psi(y, g_1(x), \dots, g_k(x)) & x \in S,\\
    \emptyset & x \notin S.
  \end{cases}
\]
Then $h$ is a set by Replacement, concentrated on $u$, so
$h \in \mathbb{F}$; and $\llbracket \psi(h, g_1, \dots, g_k) \rrbracket \supseteq S \in W$, hence
$\llbracket \psi(h, g_1, \dots, g_k) \rrbracket \in W$. By the induction hypothesis, $f(\psi)(\hat h, \hat g_1, \dots, \hat g_k)$ holds, witnessing the existential with $\hat h \in D$.
\end{proof}

\begin{corollary}[Diagonal Embedding Scheme]\label{diag}
Let $j(a) := \hat c_a$, where $c_a(x) = a$ for all $x \in \omega^A$.
Then $j$ is a $T$-provably total injective definable class function
$V \to D$, and for each formula $\varphi(x_1, \ldots, x_k)$,
\[
  T \vdash \forall a_1, \ldots, a_k\,
  \bigl( \varphi(a_1, \ldots, a_k) \gdw
   f(\varphi)(j(a_1), \ldots, j(a_k)) \bigr).
\]
In particular $f$ is a faithful self-interpretation of every $T \supseteq \ZFC$, and for every $\M \vDash T$ the map
$j^{\M} : \M \to \M_f$ is elementary (in the usual sense).
\end{corollary}

\begin{proof}
$\llbracket \varphi(c_{a_1}, \ldots, c_{a_k}) \rrbracket$ equals $\omega^A$
or $\emptyset$ depending on whether $\varphi(\bar a)$ holds or not; both lie in
$\mathcal{X}$ with least support $\emptyset$, and
$\omega^A \in W$, $\emptyset \notin W$ since
$W_\emptyset = \{\omega^\emptyset\}$. Now apply Lemma~\ref{los}.
Injectivity is the instance $\varphi := (x_1 = x_2)$. For sentences
$\varphi$ we obtain $T \vdash \varphi \gdw f(\varphi)$, hence
$T \vdash \varphi \Gdw T \vdash f(\varphi)$.
\end{proof}

\noindent
Thus $f$ is sententially equivalent to the identity interpretation; its non-triviality, established below, is purely a matter of homotopy. The key point is that a definable isomorphism is forced to be the numeral map on the natural numbers. A least-element argument inside $\mathfrak{M}$ shows that any \emph{definable} isomorphism
would have to agree with the numeral map $j$ on $\omega$, but the new element $\hat g_a$ lies outside the range of $j$. 

\begin{proof}[Proof of Theorem~\ref{non-rigid_ZFC}]
That $f$ is a self-interpretation of $T$ is immediately given by Corollary~\ref{diag}. Write $x \mathrel{\in^*} y$ for $f(x \in y)$. We will now show that $f$ produces an element of $D$ that $\M_f$ regards as a finite ordinal, i.e., an $\in^*$-element of $j(\omega)$, that does not lie in $j[\omega]:=\{j(k):k\in\omega\}$.  Any $\M$-definable isomorphism $h:  \M\to \M_f$ would have to agree with $j$ on $\omega$ and map $\omega$ onto the $\in^*$-elements of $j(\omega)$, so this element would lie in $j[\omega]$ after all.

\noindent
We show a slightly stronger statement. Let $\M\vDash T$ and suppose towards a contradiction that a formula $\eta(x, y, \bar p)$, with parameters $\bar p$ from $\M$, defines an isomorphism $h :  \M \to \M_f$. We argue in $\M$. Then $h$ is a bijection of $V$ onto $D$, and since $h$ is an isomorphism, $\varphi(\bar a)\gdw f(\varphi)(h(\bar a))$ holds for each formula $\varphi$; by Corollary~\ref{diag}, the same holds for $j$ in place of $h$.

\medskip\noindent
By Corollary~\ref{diag}, $j(\omega)$ is the element of $D$ that $\M_f$ regards as $\omega$.

\medskip\noindent
There is a non-principal ultrafilter $U$ on $\omega$, and since $|U| \leq |\mathcal{P}(\omega)| = 2^{\aleph_0}$ we have $U = U_a$ for some $a \in A$. The coordinate function $g_a(x) := x(a)$ lies in $\mathbb{F}$ with
$\|g_a\| \subseteq \{a\}$, and $\hat g_a $ is an $\in^*$-element of $j(\omega)$ not in $j[\omega]$. Indeed, since
$g_a(x) \in \omega$ for every $x \in \omega^A$ we have
$\llbracket g_a \in c_\omega \rrbracket = \omega^A \in W$, so
$\hat g_a \mathrel{\in^*} j(\omega)$ by Lemma~\ref{los}; and for each
$k \in \omega$ the set $\llbracket g_a = c_k \rrbracket$ corresponds under
$\omega^{\{a\}} \cong \omega$ to the singleton $\{k\} \notin U_a$, whence
$\llbracket g_a = c_k \rrbracket \notin W$ and $\hat g_a \neq j(k)$.

\medskip\noindent
Next, $h$ agrees with $j$ on $\omega$. The class $E:=\{k\in\omega : h(k)\neq j(k)\}$ is a set by Separation (with parameters $\bar p$). Suppose $E\neq\emptyset$, and let $k$ be its $\in$-least element. Both $h(\emptyset)$ and $j(\emptyset)$ satisfy $f(\forall z\,\neg z\in x)$, and this condition determines a unique element of $D$, since $T\vdash\exists!x\,\forall z\,\neg z\in x$ and $f$ is a self-interpretation of $T$. Hence $h(\emptyset)=j(\emptyset)$ and $k\neq 0$. Write $k=m+1$; then $h(m)=j(m)$. Let $\sigma(x,y)$ be the formula $\forall z\,(z\in y\gdw z\in x\lor z=x)$. Since $T\vdash\forall x\,\exists!y\,\sigma(x,y)$, for every $x\in D$ there is a unique $y\in D$ with $f(\sigma)(x,y)$, and both $h(m+1)$ and $j(m+1)$ are this unique $y$ for $x=h(m)=j(m)$. So $h(k)=j(k)$, contradicting $k\in E$. Hence $E=\emptyset$.

\medskip\noindent
In the same way, $h(\omega)=j(\omega)$, as both are the unique element of $D$ that $\M_f$ regards as $\omega$. Since $h$ is onto $D$, we have $\hat g_a=h(u)$ for some $u$, and since $h$ reflects membership, $h(u)\mathrel{\in^*}h(\omega)$ yields $u\in\omega$. As $E=\emptyset$, $\hat g_a=h(u)=j(u)\in j[\omega]$, a contradiction.

\medskip\noindent
Thus in no model $\M$ of $T$ does any formula $\eta(x,y,\bar p)$, with parameters from $\M$, define an isomorphism $\M \to \M_f$. By completeness, $T$ proves for each $\eta$ that for no $\bar p$ does $\eta(x,y,\bar p)$ define an isomorphism between $\id_T$ and $f$. In particular, $f$ is non-trivial and $T$ is not rigid.
\end{proof}

\noindent
If we now consider rigid models of $\ZF$ and $\ZFC$, such as the Shepherdson--Cohen model $\langle L_{\alpha_0}, \in \rangle$ where $\alpha_0 = \min \{ \alpha \in \mathrm{Ord} : L_\alpha \vDash \mathrm{ZF} \}$ or $\langle  V_\kappa, \in \rangle$ for inaccessible $\kappa$, we obtain the following corollary (the first part is assuming that there is a transitive model of $\ZF$, the second that there is an inaccessible cardinal).

\begin{corollary}\label{standardsettheory}
$\Th(\langle L_{\alpha_0}, \in \rangle)$ and $\Th(\langle V_\kappa, \in \rangle)$ for inaccessible $\kappa$ are not rigid.
\end{corollary}

\medskip\noindent
Theorem~\ref{non-rigid_ZFC} is formulated for extensions in the language $\{\in\}$. However, there can be rigid extensions in an extended language; we provide an example of this below. Nevertheless, it can be observed that the proof of Theorem~\ref{non-rigid_ZFC} holds just as well for extended languages $\mathcal{L}$, provided that the extended vocabulary also satisfies the schemata of $\ZFC$. In that case the construction interprets each new $R$ via $R_W(\hat g_1,\dots,\hat g_k) :\gdw \llbracket R(g_1,\dots,g_k) \rrbracket \in W$. This is well-defined because $\llbracket R(\bar g)\rrbracket$ is a set by $\mathcal{L}$-separation and again concentrated, and hence an element of $\mathcal{X}$. Moreover, under this hypothesis, {\L}o\'s's theorem extends to all $\mathcal{L}$-formulae, its only further scheme-dependent step being the extraction of a witnessing function in the existential case, which is an instance of $\mathcal{L}$-Collection. Thus the ultrapower satisfies every axiom of $T$ in the extended language $\mathcal{L}$ rather than only those of $T \supseteq \ZFC$ in $\{\in \}$, and the same argument as before shows the interpretation is not homotopic to the identity, so $T$ is non-rigid. 

\medskip\noindent
This observation highlights an even stronger disanalogy between $\PA$ and $\ZFC$. By Theorem~\ref{Th(N,A)rigid}, every extension $T$ of $\PA$ in a countable extended vocabulary that is true in some expansion $\fA=\langle\NN,(R_\xi)_{\xi<\alpha}\rangle$ of the standard model has a rigid extension, namely $\Th(\fA)$. However, for inaccessible $\kappa$ and arbitrary relations $(R_\xi)$ on $V_\kappa$, the expansion $\langle V_\kappa,\in,(R_\xi)\rangle$ satisfies the schemata of $\ZFC$ in the extended vocabulary, so by the observation just made its theory is not rigid.

\medskip\noindent
Since the set-theoretic counterpart to $\Th(\NN)$, the theory of hereditarily finite sets---$\Th(\langle \mathrm{HF},\in\rangle)$---is rigid (by Lemma~\ref{rigid_under_retract}), the question arises as to what threshold in axiomatic set theory precludes the possibility of a rigid extension.

\subsubsection*{Definable Enumerations}

\noindent
We can only give a partial answer to this question. First, it should be observed that the theory of a structure whose universe is enumerated by a single definable formula, and whose natural numbers are standard, is a retract of a theory covered by Theorem~\ref{Th(N,A)rigid}, and hence rigid.

\begin{definition}[Definable Enumeration]\label{def_enumeration}
Let $\mathcal{L}\supseteq\{\in\}$ be a finite language and $\mathfrak{K}$ an $\mathcal{L}$-structure. Write $\omega^{\mathfrak{K}}$ for the class of finite ordinals of $\mathfrak{K}$. We say that $\mathfrak{K}$ \emph{admits a definable enumeration} if there is a parameter-free $\mathcal{L}$-formula $\varepsilon(x,y)$ which defines in $\mathfrak{K}$ a surjective class function $e: \omega^{\mathfrak{K}}\to K$. A theory $T$ is \emph{definably enumerable} if for some such $\varepsilon$, $T$ proves that $\varepsilon$ defines a surjection of $\omega$ onto $V$.
\end{definition}

\noindent
Call $\mathfrak{K}$ \emph{$\omega$-standard} if $\omega^{\mathfrak{K}}$ is a parameter-free definable class of standard order type on which addition and multiplication are parameter-free definable.

\begin{lemma}[Enumeration Lemma]\label{enumeration_lemma}
Let $\mathcal{L}\supseteq\{\in\}$ be finite and let $\mathfrak{K}$ be an $\omega$-standard $\mathcal{L}$-structure which admits a definable enumeration. Then $\Th(\mathfrak{K})$ is rigid.\footnote{In an earlier version of this paper the argument was given only for the special case of Corollary~\ref{L_omCK_rigid} below; the general formulation was suggested by Ali Enayat and Mateusz {\L}e{\l}yk.}
\end{lemma}

\begin{proof}
Since $\mathfrak{K}$ is $\omega$-standard, $\omega^{\mathfrak{K}} \cong \N$.
Given a definable enumeration $e$, we define the least-index map
$j(x)=\min\{r\in\omega^{\mathfrak{K}} : e(r)=x\}$, which is a parameter-free
$\mathfrak{K}$-definable injection $j:  K\to\omega^{\mathfrak{K}}$. Put $D:=j[K]\subseteq \omega^{\mathfrak{K}}$ and
\[
  E(r,s):\gdw r,s\in D \;\land\; j^{-1}(r)\in^{\mathfrak{K}} j^{-1}(s),
\]
so that $D$ and $E$ are parameter-free definable over $\mathfrak{K}$ and
\begin{equation*}
  j:  \langle K,\in^{\mathfrak{K}}\rangle\;\cong\;\langle D,E\rangle .
\end{equation*}
Let $\mathfrak{A}:=\langle \NN, D,E\rangle$; if $\mathcal{L}$ contains further relation symbols $R$, add their $j$-images $R^D:=\{\langle j(a_1),\dots,j(a_k)\rangle : \langle a_1,\dots,a_k\rangle\in R^{\mathfrak{K}}\}$ to $\mathfrak{A}$ and interpret $R$ by $R^D$ below. By Theorem~\ref{Th(N,A)rigid}, $\Th(\mathfrak{A})$ is rigid.

\medskip\noindent
There is an interpretation $I:  \Th(\mathfrak{K})\to \Th(\mathfrak{A})$ whose domain-formula defines $D$ and which interprets $\in$ by $E$. The interpreted structure is isomorphic to $\mathfrak{K}$ and hence a model of $\Th(\mathfrak{K})$. Conversely, there is an interpretation $J:  \Th(\mathfrak{A}) \to \Th(\mathfrak{K})$: by $\omega$-standardness the class $\omega^{\mathfrak{K}}$ and the operations $+,\times$ are parameter-free definable in $\mathfrak{K}$, as are $D$ and $E$, and the resulting interpreted structure is precisely $\mathfrak{A}$, hence a model of $\Th(\mathfrak{A})$. The composite $J\circ I:  \Th(\mathfrak{K})\to \Th(\mathfrak{K})$ yields the definable copy $\langle D,E\rangle$ of $\mathfrak{K}$, and $j$ is a parameter-free $\mathfrak{K}$-definable isomorphism onto it; the sentence expressing that $j$ defines an isomorphism between $\id$ and $J\circ I$ is true in $\mathfrak{K}$ and hence provable in $\Th(\mathfrak{K})$. Thus $\Th(\mathfrak{K})$ is a retract of $\Th(\mathfrak{A})$, and $\Th(\mathfrak{K})$ is rigid by Lemma~\ref{rigid_under_retract}.
\end{proof}

\noindent
Lemma~\ref{enumeration_lemma} also shows that the restriction to the language $\{\in\}$ in Theorem~\ref{non-rigid_ZFC} cannot be dropped. Assume there is a transitive model of $\ZF$; then $\langle L_{\alpha_0},\in\rangle$ is a countable transitive model of $\ZFC$, so there is a bijection $R\subseteq\omega\times L_{\alpha_0}$ between $\omega$ and $L_{\alpha_0}$. The structure $\langle L_{\alpha_0},\in,R\rangle$ is $\omega$-standard and admits the definable enumeration given by $R$, so the theory $\Th(\langle L_{\alpha_0},\in, R\rangle)$ is a consistent rigid extension of $\ZFC$ in the language $\{\in,R\}$, although $\Th(\langle L_{\alpha_0},\in\rangle)$ is not rigid. This does not contradict Theorem~\ref{non-rigid_ZFC} for extended languages, as the Separation scheme for formulae mentioning $R$ fails in $\langle L_{\alpha_0},\in,R\rangle$.

\medskip\noindent
The rigidity of $\Th(\langle\mathrm{HF},\in\rangle)$ noted above is the simplest instance of Lemma~\ref{enumeration_lemma}, as Ackermann coding provides a parameter-free definable bijection between $\mathrm{HF}$ and its finite ordinals. Another instance is the theory of the rigid model $\langle \Lom, \in \rangle$, where $\Lom$ denotes the least admissible set containing $\omega$. $\langle \Lom, \in \rangle$ satisfies $\KP$, $V=L$, and that every element is countable.\footnote{Following \citet{Mathias2001} and \citet{EnayatLelyk2024}, $\KP$ denotes Kripke--Platek set theory without Infinity and with the Foundation scheme restricted to $\Pi_1$-formulae, equivalently with the scheme of $\in$-induction restricted to $\Sigma_1$-formulae. Where the full scheme is needed, we add it as ``Full Set Induction'' (equivalently, the full Foundation scheme). In this formulation $\KP$ is finitely axiomatizable.}

\begin{corollary}\label{L_omCK_rigid}
$\Th(\langle \Lom, \in \rangle)$ is rigid.
\end{corollary}

\begin{proof}
Over $\Lom$ there is a definable injection $i :  \Lom \to \omega$; equivalently, a definable bijection of $\Lom$ onto $\mathrm{ran}(i) \subseteq\omega$ \citep[cf.][Ch.~5, Cor.~5.5]{Barwise2017}. So $\langle\Lom,\in\rangle$ admits a definable enumeration. Moreover, $\omega$ is a definable element of $\Lom$ and $0,1,+,\times$ are $\Delta_1$ over $\langle\Lom,\in\rangle$ by $\Sigma$-recursion; as $\Lom$ is transitive, $\langle\Lom,\in\rangle$ is $\omega$-standard. The claim follows by Lemma~\ref{enumeration_lemma}.
\end{proof}

\noindent
The same holds for every admissible $\alpha > \omega$ such that $\langle L_\alpha, \in\rangle$ admits a parameter-free definable enumeration, e.g., whenever every element of $L_\alpha$ is $\Sigma_1$-definable in $L_\alpha$ without
parameters. Besides $\omCK$, this includes the least stable ordinal $\sigma_0$, since $L_{\sigma_0}$ consists exactly of the parameter-free $\Sigma_1$-definable
elements of $L$ \citep[Ch.~5, Thm.~7.10]{Barwise2017}. Note that this condition is stronger than projectibility in the sense of Barwise\footnote{Cf.\ \cite[p.~188]{Barwise2017}.}. Moreover, rigidity is not monotone along the $L$-hierarchy: $\alpha_0 < \sigma_0$, but $\Th(\langle L_{\alpha_0}, \in\rangle)$ is not rigid.

\subsubsection*{The Cantor Threshold}

\medskip\noindent
The preceding subsection motivates the study of the class of consistent theories $T\supseteq \KP+\mathrm{Infinity}$ in $\{\in\}$ such that for every formula $\varepsilon(x,y)$, $T$ proves that $\varepsilon$ does not define a surjection of $\omega$ onto $V$. In particular such a theory is not definably enumerable, as none of its models admit a definable enumeration. We might say that such a theory $T$ has \emph{crossed the Cantor threshold}. A sufficient condition for this is $T$ proving the full Separation scheme or the full Collection scheme. Indeed, let $\varepsilon(x,y)$ define a total class function $e$ on $\omega$ in a model of $T$. By Separation, the set $d := \{n \in \omega : n \notin e(n)\}$ exists; if $e$ were onto, then $d = e(k)$ for some $k$, and $k \in d \leftrightarrow k \notin d$. Alternatively, Collection applied to $\omega$ and $\varepsilon$ yields a set $b$ containing every $e(n)$; if $e$ were onto, then $b \in b$, contradicting Foundation. For theories such as $\Ze$, $\ZF^-$, and $\ZF$ the method of Lemma~\ref{enumeration_lemma} is therefore unavailable. Nevertheless, every bounded fragment of $\ZFC$ still has rigid extensions:\footnote{Theorem~\ref{EL_rigid_fragments} is due jointly to Ali Enayat, Piotr Gruza, and Mateusz {\L}e{\l}yk; Gruza arrived independently at the idea of obtaining rigid theories from countably coded models, that is, from $\omega$-standard models admitting a definable enumeration. The proof given here follows a proof by Enayat for $\Th_{\Pi_{n+1}}(\ZF)$.  We are much indebted to Enayat and {\L}e{\l}yk for pointing out an error in an earlier version of this section, and for allowing us to include the theorem and its proof here.}

\begin{theorem}[Enayat, Gruza, {\L}e{\l}yk]\label{EL_rigid_fragments}
Assume that $\ZF$ has a well-founded model. Then for every $n\in\omega$ there is a complete rigid theory $T_n$ such that
\[
  T_n\;\supseteq\;\KP+\text{Full Set Induction}+\Th_{\Pi_{n+1}}(\ZF+V=L),
\]
where $\Th_{\Pi_{n+1}}(\ZF+V=L)$ is the set of $\Pi_{n+1}$-sentences provable in $\ZF+V=L$. In particular, $T_n\supseteq\Th_{\Pi_{n+1}}(\ZFC)$.
\end{theorem}

\begin{proof}[Proof sketch]
Fix a well-founded $\M\vDash\ZF+V=L$ (pass to its constructible universe if necessary) and $n\in \omega$, let $m:=\max\{3,n\}$ and let $\mathfrak{K}:=\mathfrak{K}_m(\M)$ be the submodel of $\M$ whose universe consists of the parameter-free $\Sigma_m$-definable elements of $\M$. By \citet[Theorem~3(a)]{EnayatLelyk2026corr}, $\mathfrak{K}$ is a $\Sigma_m$-elementary (equivalently, $\Pi_m$-elementary) substructure of $\M$. Hence every
$\Pi_{m+1}$-sentence true in $\M$ is true in $\mathfrak{K}$; in particular,
$\mathfrak{K}\vDash\Th_{\Pi_{m+1}}(\ZF+V=L)$. Since $n\leq m$, this gives
$\mathfrak{K}\vDash\Th_{\Pi_{n+1}}(\ZF+V=L)$. Moreover, since $m\geq 3$, every $\Pi_3$-sentence true in $\M$ is true in $\mathfrak{K}$; as $\KP$ is axiomatized by a set of $\Pi_3$-sentences, all of which are provable in $\ZF$, we also obtain $\mathfrak{K}\vDash\KP$. Since $\mathfrak{K}$ is a substructure of the well-founded $\M$ with the inherited membership relation, $\in^{\mathfrak{K}}$ is externally well-founded; therefore $\mathfrak{K}$ satisfies the full scheme of set induction and $\omega^{\mathfrak{K}}$ has standard order type. Since moreover $\mathfrak{K}\vDash\KP+\text{Infinity}$, $\omega$ is $\Delta_0$-definable in $\mathfrak{K}$ and $+,\times$ are definable by $\Sigma$-recursion, so $\mathfrak{K}$ is $\omega$-standard. By \citet[Theorem~3(b)]{EnayatLelyk2026corr}, $\mathfrak{K}_m(\mathfrak{K})=\mathfrak{K}$, i.e., every element of $\mathfrak{K}$ is parameter-free $\Sigma_m$-definable \emph{in $\mathfrak{K}$}. Using the $\Sigma_m$-satisfaction predicate available in $\mathfrak{K}$, the map sending (the code of) a $\Sigma_m$-formula in one free variable to the unique element it defines, and to $\emptyset$ otherwise, is a parameter-free definable class function $e: \omega\to K$. So $\mathfrak{K}$ admits a definable enumeration, and $T_n:=\Th(\mathfrak{K})$ is rigid by Lemma~\ref{enumeration_lemma}.
\end{proof}

\noindent
The picture for axiomatic set theory that now emerges is the following. The familiar axiomatized fragments of $\ZFC$ are not themselves rigid. If $\ZFC$ is consistent, no finitely axiomatized subtheory of $\ZFC$ is rigid, since by Lemma~\ref{upwards_preservation_rigidity} its rigidity would entail that of $\ZFC$, contradicting Theorem~\ref{non-rigid_ZFC}; this applies, for instance, to $\KP$. More generally, no finitely axiomatized theory $S$ in the language $\{\in\}$ that is consistent with $\ZFC$ is rigid, since otherwise $S+\ZFC$ would be rigid by Lemma~\ref{upwards_preservation_rigidity}. Likewise, no consistent r.e.\ subtheory of $\ZFC$ that is inductive sequential, such as $\Ze$, $\ZF^-$, or $\KP+\text{Full Set Induction}$, is rigid for the reasons mentioned in footnote~\ref{Enayat_footnote}. However, every bounded fragment $\Th_{\Pi_n}(\ZFC)$ of $\ZFC$ has a consistent rigid extension, assuming that $\ZF$ has a well-founded model. In particular, for any fixed $m$, $\Sigma_m$-Separation and $\Sigma_m$-Collection can be added (each instance of these schemata is a $\Pi_{m+3}$-sentence provable in $\ZF$, so it suffices to take $T_n$ with $n\geq m+2$). For weak fragments such as $\KP$, rigid extensions exist even without this assumption (Corollary~\ref{L_omCK_rigid}). Where between the bounded fragments and $\ZFC$ rigid extensions disappear remains open. In particular, we do not know whether crossing the Cantor threshold already excludes them, nor whether rigidity in set theory can be obtained by any means other than definable enumerations. We therefore close with a question that we find of philosophical as well as mathematical interest.

\medskip\noindent
\textbf{Open Question.}
\begin{enumerate}
  \item Is there a consistent rigid theory that has crossed the Cantor
        threshold?
  \item In particular, is there a consistent rigid extension of $\Ze$,
        of $\ZF^-$, or of $\ZF$ in the language $\{\in\}$?
        (Any rigid
        extension of $\ZF$ must refute $\mathrm{AC}$.)
\end{enumerate}

\section{Conclusion}

The results of this work show that rigidity of theories is a distinct and interesting invariance property that cannot be reduced in any direct manner to the invariance properties already well understood in model theory and proof theory.

\medskip\noindent
There is no direct overlap with categoricity in a cardinality, nor with rigidity of intended models, nor with completeness of theories.
Our selection of examples has already shown that we do indeed have $\aleph_1$-categorical theories that are rigid ($\ACF_0$), but also rigid theories that are not categorical in any cardinality ($\RCF$). Furthermore, while we have theories of rigid models that are also rigid ($\RCF$, True Arithmetic), we also have rigid theories that do not possess rigid models ($\ACF$), as well as theories of rigid models that are not rigid (Corollary~\ref{standardsettheory}). Moreover, our investigations reveal that there are incomplete theories which are rigid ($\ACF$).   

\medskip\noindent
What rigidity of theories captures—and what none of the other properties mentioned capture—is the extent to which the axioms of a theory constrain the interpretation of its own primitive symbols under the transformations inherent to the category of theories, i.e., how much freedom the expressive and inferential resources of the theory give one when reinterpreting the primitives of the theory. 

\medskip\noindent
Compared to more closely related concepts of semantic determinacy, such as solidity,\footnote{Cf.\ \cite{EnayatLelyk2024}.} it can be observed that rigidity also allows for finer distinctions. Both $\PA$ and $\RCF$ are solid, but only the latter is rigid. Note that rigidity of a \emph{complete} or \emph{finitely axiomatized} theory implies solidity (without parameters); solidity (without parameters), in turn, implies that a theory (in a finite language) is rigid with respect to retractions. The possibility of making finer distinctions becomes even more apparent when we consider the existence of rigid extensions. Both $\PA$ and $\ZFC$ are solid and thus rigid with respect to (co-)retractions, but only the latter lacks a consistent, rigid extension in its own language, whereas $\Th(\NN)$ is rigid. It is not unlikely that these classifications also correspond to an intuitive assessment of a certain kind of semantic determinacy among the theories mentioned, which could be described as \emph{conceptual determinacy}; however, this requires further analysis.

\medskip\noindent
This brief overview of examples already demonstrates that rigidity of theories is a multifaceted and useful property that warrants further investigation. A more in-depth philosophical examination of this concept and a detailed systematic study of the general number of self-interpretations of theories will, as already mentioned, be the subject of a separate article.
 
\begingroup
\small
\setstretch{1}
\setlength{\parindent}{0pt}
\setlength{\parskip}{0.4em}

\textbf{Funding Acknowledgements.}
Funded by the European Union (ERC, FORMALISM, 101044114).
Views and opinions expressed are, however, those of the author(s)
only and do not necessarily reflect those of the European Union
or the European Research Council Executive Agency. Neither the
European Union nor the granting authority can be held responsible
for them.

\textbf{AI Statement.}
The authors used Claude Opus and ChatGPT Astra to help find relevant
literature on the topic and to assist with developing proofs and
checking their accuracy, drafting parts of text, and identifying typos.
All mathematical content suggested by LLMs was independently
checked by the authors. The authors take full responsibility
for the content of the paper.
\par
\endgroup


\begin{thebibliography}{35}
\providecommand{\natexlab}[1]{#1}
\providecommand{\url}[1]{\texttt{#1}}
\expandafter\ifx\csname urlstyle\endcsname\relax
  \providecommand{\doi}[1]{doi: #1}\else
  \providecommand{\doi}{doi: \begingroup \urlstyle{rm}\Url}\fi

\bibitem[Ajdukiewicz(1934)]{ajdukiewicz1934}
Kazimierz Ajdukiewicz.
\newblock {Sprache und Sinn}.
\newblock \emph{Erkenntnis}, 4\penalty0 (1):\penalty0 100--138, 1934.

\bibitem[Apelt(1966)]{apelt1966}
Harry Apelt.
\newblock Axiomatische {U}ntersuchungen {\"u}ber einige mit der
  {P}resburgerschen {A}rithmetik verwandte {S}ysteme.
\newblock \emph{Mathematical Logic Quarterly}, 12\penalty0 (1):\penalty0
  131--168, 1966.

\bibitem[Barwise(2017)]{Barwise2017}
Jon Barwise.
\newblock \emph{Admissible Sets and Structures}, volume~7 of
  \emph{Perspectives in Logic}.
\newblock Cambridge University Press, Cambridge, 2017.
\newblock Reissue of the 1975 Springer edition.

\bibitem[Bays and Peterzil(2019)]{bays2019definability}
Martin Bays and Ya'acov Peterzil.
\newblock Definability in the group of infinitesimals of a compact Lie group.
\newblock \emph{Confluentes Mathematici}, 11(2):3--23, 2019.

\bibitem[B{\`e}s(2008)]{bes2008application}
Alexis B{\`e}s.
\newblock An application of the {Feferman--Vaught} theorem to automata and
  logics for words over an infinite alphabet.
\newblock \emph{Logical Methods in Computer Science}, 4, 2008.

\bibitem[Boolos et~al.(2007)Boolos, Burgess, and Jeffrey]{boolos2007}
George~S. Boolos, John~P. Burgess, and Richard~C. Jeffrey.
\newblock \emph{Computability and Logic}.
\newblock Cambridge University Press, New York, 5th edition, 2007.

\bibitem[Carnap(1937)]{carnap1937}
Rudolf Carnap.
\newblock \emph{The Logical Syntax of Language}.
\newblock Kegan Paul, Trench, Trubner \& Co., London, 1937.

\bibitem[Cegielski(1981)]{cegielski1981}
Patrick Cegielski.
\newblock Th{\'e}orie {\'e}l{\'e}mentaire de la multiplication des entiers
  naturels.
\newblock In \emph{Model Theory and Arithmetic: Comptes Rendus d'une Action
  Th{\'e}matique Programm{\'e}e du CNRS sur la Th{\'e}orie des Mod{\`e}les et
  l'Arithm{\'e}tique, Paris, France, 1979/80}, pages 44--89. Springer, 1981.


\bibitem[Chang and Keisler(1990)]{ChangKeisler1990}
Chen~Chung Chang and H.~Jerome Keisler.
\newblock \emph{Model Theory}, volume~73 of \emph{Studies in Logic and the
  Foundations of Mathematics}.
\newblock North-Holland, Amsterdam, 3rd edition, 1990.

\bibitem[Enayat(2016)]{Enayat2016}
Ali Enayat.
\newblock Variations on a {V}isserian theme.
\newblock In Jan van Eijck, Rosalie Iemhoff, and Joost~J. Joosten, editors,
  \emph{Liber {A}micorum {A}lberti: a tribute to {A}lbert {V}isser}, pages
  99--110. College Publications, London, 2016.

\bibitem[Enayat et~al.(2011)Enayat, Schmerl, and Visser]{Enayat2011}
Ali Enayat, James~H. Schmerl, and Albert Visser.
\newblock $\omega$-models of finite set theory.
\newblock In Juliette Kennedy and Roman Kossak, editors, \emph{Set Theory,
  Arithmetic, and Foundations of Mathematics: Theorems, Philosophies},
  volume~36 of \emph{Lecture Notes in Logic}, pages 43--65. Cambridge
  University Press, Cambridge, 2011.

\bibitem[Enayat and {\L}e{\l}yk(2024)]{EnayatLelyk2024}
Ali Enayat and Mateusz {\L}e{\l}yk.
\newblock Categoricity-like properties in the first order realm.
\newblock \emph{Journal for the Philosophy of Mathematics}, 1:63--98, 2024.

\bibitem[Enayat and {\L}e{\l}yk(2026)]{EnayatLelyk2026corr}
Ali Enayat and Mateusz {\L}e{\l}yk.
\newblock Corrigendum \& addendum to ``{C}ategoricity-like properties in the
  first order realm''.
\newblock arXiv preprint arXiv:2603.06711, 2026.


\bibitem[Feferman(1958)]{Feferman1958}
Solomon Feferman.
\newblock Arithmetically definable models of formalized arithmetic.
\newblock \emph{Notices of the American Mathematical Society}, 5:\penalty0 679,
  1958.

\bibitem[Feferman and Vaught(1959)]{feferman1959first}
Solomon Feferman and Robert Vaught.
\newblock The first order properties of products of algebraic systems.
\newblock \emph{Fundamenta Mathematicae}, 47(1):57--103, 1959.

\bibitem[Hedrl{\'\i}n and Pultr(1966)]{hedrlin1966rigid}
Zden{\v e}k Hedrl{\'\i}n and Ale{\v s} Pultr.
\newblock On rigid undirected graphs.
\newblock \emph{Canadian Journal of Mathematics}, 18:\penalty0 1237--1242,
  1966.


\bibitem[Hilbert(1899)]{hilbert1899}
David Hilbert.
\newblock \emph{Grundlagen der Geometrie}.
\newblock Festschrift zur Feier der Enthüllung des Gauss-Weber-Denkmals in
  Göttingen, B.G. Teubner, Leipzig, 1899.

\bibitem[Hinman(2005)]{hinman2005}
Peter~G. Hinman.
\newblock \emph{Fundamentals of Mathematical Logic}.
\newblock A K Peters, Wellesley, Massachusetts, 2005.

\bibitem[Hodges(1993)]{hodges1993}
Wilfrid Hodges.
\newblock \emph{Model Theory}.
\newblock Cambridge University Press, Cambridge, England, 1993.

\bibitem[Jahnke and Kartas(2025)]{jahnke2025beyond}
Franziska Jahnke and Konstantinos Kartas.
\newblock Beyond the {F}ontaine-{W}intenberger theorem.
\newblock \emph{Journal of the American Mathematical Society},
  38(4):997--1047, 2025.
  
\bibitem[Jeř{\'a}bek(2012)]{Jerabek2012}
Emil Jeř{\'a}bek.
\newblock Sequence encoding without induction.
\newblock \emph{Mathematical Logic Quarterly}, 58\penalty0 (3):\penalty0
  244--248, 2012.

\bibitem[Kamio and Hora(2026)]{kamio2026solution}
Yuhi Kamio and Ryuya Hora.
\newblock Solution to {L}awvere's first problem: A {G}rothendieck topos that
  has proper class many quotient topoi.
\newblock \emph{Advances in Mathematics}, 487:\penalty0 110751, 2026.

\bibitem[Kanovei and Shelah(2004)]{KanoveiShelah2004}
Vladimir Kanovei and Saharon Shelah.
\newblock A definable nonstandard model of the reals.
\newblock \emph{The Journal of Symbolic Logic}, 69\penalty0 (1):\penalty0
  159--164, 2004.

\bibitem[Kaye(1991)]{Kaye1991}
Richard Kaye.
\newblock \emph{Models of Peano Arithmetic}, volume~15 of \emph{Oxford Logic
  Guides}.
\newblock Oxford University Press, Oxford, 1991.

\bibitem[Ketelsen(2026)]{ketelsen2026understanding}
Margarete Ketelsen.
\newblock \emph{Understanding the model theory of henselian valued fields
  by parts}.
\newblock PhD thesis, University of M\"unster, 2026.

\bibitem[Marker(2000)]{marker2000}
David Marker.
\newblock Introduction to model theory.
\newblock In Deirdre Haskell, Anand Pillay, and Charles Steinhorn, editors,
  \emph{Model Theory, Algebra, and Geometry}, volume~39 of
  \emph{Mathematical Sciences Research Institute Publications},
  pages 15--35.
\newblock Cambridge University Press, Cambridge, 2000.

\bibitem[Mathias(2001)]{Mathias2001}
A.~R.~D. Mathias.
\newblock The strength of {M}ac {L}ane set theory.
\newblock \emph{Annals of Pure and Applied Logic}, 110\penalty0 (1--3):\penalty0
  107--234, 2001.

\bibitem[Milne(2022)]{milneFT}
James S. Milne.
\newblock \emph{Fields and Galois Theory}, version~5.10.
\newblock Course notes, 2022.
\newblock Available at \url{https://www.jmilne.org/math/CourseNotes/FT.pdf}.

\bibitem[Mostowski(1952)]{mostowski1952direct}
Andrzej Mostowski.
\newblock On direct products of theories.
\newblock \emph{The Journal of Symbolic Logic}, 17\penalty0 (1):\penalty0
  1--31, 1952.

\bibitem[Orey(1961)]{Orey1961}
Steven Orey.
\newblock Relative interpretations.
\newblock \emph{Zeitschrift f{\"u}r mathematische Logik und Grundlagen der
  Mathematik}, 7:\penalty0 146--153, 1961.

\bibitem[Otero et~al.(1996)Otero, Peterzil, and Pillay]{Pillay1996}
Margarita Otero, Ya'acov Peterzil, and Anand Pillay.
\newblock On groups and rings definable in o–minimal expansions of real
  closed fields.
\newblock \emph{Bulletin of the London Mathematical Society}, 28\penalty0
  (1):\penalty0 7--14, 1996.

\bibitem[Pakhomov and Zapryagaev(2020)]{Pakhomov2020}
Fedor Pakhomov and Alexander Zapryagaev.
\newblock Multi-dimensional interpretations of {P}resburger arithmetic in itself.
\newblock \emph{Journal of Logic and Computation}, 30\penalty0 (8):\penalty0
  1681--1693, 2020.
  
\bibitem[Poizat(1983)]{poizat1983theorie}
Bruno Poizat.
\newblock Une th\'eorie de Galois imaginaire.
\newblock \emph{The Journal of Symbolic Logic}, 48(4):1151--1170, 1983.

\bibitem[Poizat(1988)]{poizat1988mm}
Bruno Poizat.
\newblock MM. Borel, Tits, Zil'ber et le G\'en\'eral Nonsense.
\newblock \emph{The Journal of Symbolic Logic}, 53(1):124--131, 1988.

\bibitem[Poizat(2026)]{poizat2026parametres}
Bruno Poizat.
\newblock Param\`etres dans les corps alg\'ebriquement clos: Une
  ex\'eg\`ese mod\`ele-th\'eorique d'un r\'esultat de Monsieur Aleksandr
  Vasil'evi\'c Borovik, suivie d'un commentaire de m\^eme nature sur le
  Th\'eor\`eme de Borel-Tits.
\newblock \emph{Annales math\'ematiques du Qu\'ebec},
  50(1):295--321, 2026.

\bibitem[Schweikardt(2005)]{schweikardt2005arithmetic}
Nicole Schweikardt.
\newblock Arithmetic, first-order logic, and counting quantifiers.
\newblock \emph{ACM Transactions on Computational Logic (TOCL)}, 6\penalty0
  (3):\penalty0 634--671, 2005.

\bibitem[Scott(1960)]{Scott1960}
Dana Scott.
\newblock On a theorem of {Rabin}.
\newblock \emph{Indagationes Mathematicae (Proceedings)}, 63:\penalty0 481--484, 1960.

\bibitem[Smory{\'n}ski(2012)]{smorynski2012logical}
Craig Smory{\'n}ski.
\newblock \emph{Logical number theory I: An introduction}.
\newblock Springer Science \& Business Media, 2012.

\bibitem[Stonestrom(2022)]{stonestrom2022some}
Atticus Stonestrom.
\newblock Some model theory of {$\Th(\mathfrak{N}, \times)$}.
\newblock \emph{Mathematical Logic Quarterly}, 68\penalty0 (3):\penalty0
  288--303, 2022.

\bibitem[Tarski et~al.(1953)Tarski, Mostowski, and Robinson]{tarski1953}
Alfred Tarski, Andrzej Mostowski, and Raphael~M Robinson.
\newblock \emph{Undecidable Theories}.
\newblock North-Holland Publishing Company, Amsterdam, 1953.


\bibitem[Visser(2005)]{visser2005faith}
Albert Visser.
\newblock Faith \& falsity.
\newblock \emph{Annals of Pure and Applied Logic}, 131(1--3):103--131, 2005.


\bibitem[Visser(2006)]{Visser2006}
Albert Visser.
\newblock Categories of theories and interpretations.
\newblock In Ali Enayat, Iraj Kalantari, and Mojtaba Moniri, editors,
  \emph{Logic in Tehran}, volume~26 of \emph{Lecture Notes in Logic}, pages
  284--341. Association for Symbolic Logic, 2006.
\newblock \doi{10.1017/9781316755747.019}.

\bibitem[Visser(2017{\natexlab{a}})]{Visser2017}
Albert Visser.
\newblock On {Q}.
\newblock \emph{Soft Computing}, 21\penalty0 (1):\penalty0 39--56, 2017.

\bibitem[Visser(2017{\natexlab{b}})]{visser2017interpretation}
Albert Visser.
\newblock The interpretation existence lemma.
\newblock In Gerhard J{\"a}ger and Wilfried Sieg, editors, \emph{Feferman on
  Foundations: Logic, Mathematics, Philosophy}, volume~13 of \emph{Outstanding
  Contributions to Logic}, pages 101--144. Springer, Cham, 2017.

\bibitem[Visser(2019)]{Visser2019}
Albert Visser.
\newblock The small-is-very-small principle.
\newblock \emph{Mathematical Logic Quarterly}, 65\penalty0 (4):\penalty0
  453--478, 2019.

\bibitem[Zapryagaev(2023)]{ZapThesis}
Alexander~Alexandrovich Zapryagaev.
\newblock \emph{Interpretations in weak arithmetical theories}.
\newblock PhD thesis, National Research University Higher School of Economics,
  Moscow, 2023.

\bibitem[Zoethout(2015)]{zoethout2015interpretations}
Jetze Zoethout.
\newblock \emph{Interpretations in {P}resburger {A}rithmetic}.
\newblock Bachelor's thesis, Utrecht University, Utrecht, 2015.

\end{thebibliography}
\end{document}